# Minimum Sensor Placement for Single Robust Observability of Structured Complex Networks

Xiaofei Liu[1,*]     Sérgio Pequito[1,2,*]     Soummya Kar[1]

Bruno Sinopoli[1]     A. Pedro Aguiar[2,3]

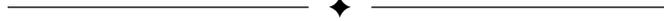

### Abstract

This paper addresses problems on the robust structural design of complex networks. More precisely, we address the problem of deploying the minimum number of dedicated sensors, i.e., those measuring a single state variable, that ensure the network to be structurally observable under disruptive scenarios. The disruptive scenarios considered are as follows: (i) the malfunction/loss of one arbitrary sensor, and (ii) the failure of connection (either unidirectional or bidirectional communication) between a pair of agents. First, we show these problems to be NP-hard, which implies that efficient algorithms to determine a solution are unlikely to exist. Secondly, we propose an intuitive two step approach: (1) we achieve an arbitrary minimum sensor placement ensuring structural observability; (2) we develop a sequential process to find minimum number of additional sensors required for robust observability. This step can be solved by recasting it as a weighted set covering problem. Although this is known to be an NP-hard problem, feasible approximations can be determined in polynomial-time, and used to obtain feasible approximations to the robust structural design problems with optimality guarantees. Finally, we discuss how the proposed methodology can be extended to the multiple sensor/link failure.

*The first two authors made equal contributions to the research and writing of this paper.

The work was partially supported by grant SFRH/BD/33779/2009, from Fundação para a Ciência e a Tecnologia (FCT) and the CMU-Portugal (ICTI) program.

[1] Department of Electrical and Computer Engineering, Carnegie Mellon University, Pittsburgh, PA 15213

[2] Institute for Systems and Robotics, Instituto Superior Técnico, Technical University of Lisbon, Lisbon, Portugal

[3] Research Center for Systems and Technologies, Faculty of Engineering, University of Porto, Porto, Portugal

Emails: {xiaofeil,spequito,soummyak,brunos}@cmu.edu, and pedro.aguiar@fe.up.pt

# 1 INTRODUCTION

State observability is an important prerequisite for the operation of complex networks with particular emphasis on critical infrastructures such as electric power grid, manufacturing plants, and transportation systems, just to name a few. Over the past decade, the complexity of these infrastructures has grown, accompanied by the increasing likelihood of failures that will disrupt the normal operation of the system. These faults can be the result of an external malicious party that aims to disrupt the operation of the system [1], [2], [3], or occur due to natural causes, such as malfunction of systems components [4].

Consequently, it is of major importance to enforce the *resilience* of these systems by proper deployment of sensors in the network. Hereafter, by exploring the structural *vulnerabilities* of the network, we address the problem of placing the minimum number of sensors that not only ensure observability [5], [6], [7], [8], but also ensure such properties under disruptive scenarios:

**$P_s$**  malfunction/loss of one arbitrary sensor;

**$P_l$**  failure of connection (either unidirectional or bidirectional) between a pair of subsystems (also referred to as *agents*) in an interconnected dynamical system (also referred to as a *complex network*).

Note that in this paper, an unidirectional connection is also referred to as a *directed link*, and a bidirectional connection is also referred to as an *undirected link*.

The networks considered hereafter are described by a (possibly large) dynamical linear time-invariant (LTI) system given by

$$\dot{x}(t) = Ax(t), \quad x(0) = x_0 \in \mathbb{R}^n,  \tag{1}$$

where $x \in \mathbb{R}^n$ denotes the state of the system, and $A \in \mathbb{R}^{n \times n}$ is the dynamic matrix. Hereafter, we concentrate our attention in determining the minimum placement of *dedicated* sensors, i.e., each sensor measures a single state variable, such that the observability of the system is attainable under the disruptive scenarios $P_s$ and $P_l$. The dedicated sensor placement is described by

$$y(t) = \mathbb{I}_n(\mathcal{J})x(t),  \tag{2}$$

where $y \in \mathbb{R}^q$ denotes the measured outputs of the system, $\mathbb{I}_n(\mathcal{J})$ corresponds to the $q \times n$ matrix comprising by all $j$th rows of the $n \times n$ identity matrix $\mathbb{I}_n$ with $j \in \mathcal{J} = \{\tau_1, \ldots, \tau_q\}$, where $\tau_i$ is the $i$th measured state variable. Dedicated sensors are common in different complex networks: (i) in the electric power grid the measurements can consist of the frequency of a bus or the power consumed by an aggregate load; and (ii) in multi-agent networks such sensors may measure the state of an agent.

*Related Work*

First, we notice that due to the duality between controllability and observability in LTI systems, the results about the former can be restated in terms of the latter. The characterization of the minimum collection of state variables that need to be measured to attain observability, also known as *minimum observability problem*, was studied in [5] and [9], where in the former the NP-hardness of the problem was established and latter exact solutions could be determined to systems whose dynamics was characterized by a simple matrix. Under the additional Grammian-based energy constraint, in [10] the authors showed that an important class of metrics had properties that allows efficient global optimization for sensor placement, while providing some optimality guarantees and when the initial placement ensures observability. Later, in [11] this framework was extended to the case where an initial placement was not required, and some additional energy constraints were considered.

Alternatively, because more often than not the parametrization of the system's dynamics is not accurately known, a natural direction is to consider structural systems [12]. These allow to consider only the fixed zero or independent parameter patterns of the system plant matrices, and to establish control theoretical properties such as *structural observability*. A system is structurally observable if for almost all parameterizations, satisfying a given pattern, the system is observable in the classical sense [12]. In [6], the structural minimum observability problem was shown to be polynomially solvable, and the characterization of all feasible solutions was provided. More recently, in [13] the structural minimum observability problem was considered with heterogeneous measuring costs and also shown to be polynomially solvable; in addition, in [7] computational methods with lower complexity were provided to the case where the cost was binary.

Although, if the minimum structural observability problem is constrained to an initial collection of sensors (possibly measuring more than a single state variable) as initially studied in [14], then the problem becomes NP-hard [15].

Several verification conditions for $\mathbf{P_s}$ and $\mathbf{P_1}$ have been previously proposed in the literature. For instance, the robustness with respect to sensor failure was initially discussed in [16], and in [17] (and references therein) the problem of recovering structural observability by allocating additional sensing capabilities was considered. Subsequently, the implications of the later into fault detection schemes was explored in [18]. In addition, our work also contrasts with those addressing robustness with respect to link failure. For instance, in [19], [16] the authors assessed the impact of directed edges failure in the structural observability of the system. More specifically, in [19], the authors provided graph-theoretic procedures to identify the minimal sets of links which were essential for preserving the structural observability. On the other hand, in [16] generic conditions were provided for the case where the link failures concerned the interconnection between subsystems. More recently, in [20], [21], the relationship between robust observability and network features such as topological transitivity and degree of the network was studied, respectively. For systems with leader-follower type architecture, where the dynamic matrix has non-zero diagonal entries, [22] introduced the notions of agent and link observability indices to characterize the importance of individual links on preserving the observability of the overall network.

These results provide verification conditions and procedures to determine robustness with respect to link failures. In this paper, we extend preliminary results presented in [23], where the dynamic matrix was assumed to be irreducible and symmetric, to the case of arbitrary dynamic matrices. Notice that it is often the case that the dynamic matrix is not symmetric nor irreducible; e.g., power systems [24], directed multi-agent networks [25], and other large-scale systems [26]. Hence, justifying the extension in the current manuscript for more general scenarios where the design for robustness/resilience is essential. In addition, we notice that from a technical perspective the methods employed are more elaborate than those in [23], due to the increased complexity of the problem. Further, we have also included a discussion on how the

                                                                                    

proposed methods can be used to tackle the problem when multiple sensor/link failures occur simultaneously. Finally, we notice that a natural robust feasible (generally suboptimal) solution consists of taking the union of two minimal and disjoint collections of sensors ensuring structural observability, which could be obtained by using the algorithms in [13] for determining minimal dedicated sensor configurations achieving structural observability under heterogeneous sensor cost constraints: the first collection is obtained considering uniform sensor cost, and the second by imposing infinite cost to the variables in the first collection. Nevertheless, the solution is not generally optimal, and can be loose when sparse networks are considered, as it will become clear in the sequel. ○

The main contributions of this paper are as follows: we address the problem of placing the minimum number of dedicated sensors that ensure the network to be structurally observable under disruptive scenarios. The two disruptive scenarios considered are as follows: (i) the malfunction/loss of one arbitrary sensor, and (ii) the failure of link (either directed or undirected) between agents in a complex network. We show that both problems are NP-hard, which implies that they are unlikely to be polynomially solvable. Therefore, we propose an intuitive two step approach: first, we achieve an arbitrary minimum sensor placement configuration ensuring structural observability; secondly, we develop a sequential process to obtain the minimum number of additional sensors required for robust observability with respect to any failure. This step can be solved by recasting it as a weighted set covering problem. Although this is known to be an NP-hard problem, feasible approximations can be determined in polynomial-time that can be used to obtain a feasible approximation to the robust structural design problems with optimality guarantees. Nevertheless, we show that the sensor placement considered in the first step may influence the size of the optimal result that we achieve in the second step. Subsequently, designing the system with respect to an arbitrary number of sensor/link failures is at least as computationally difficult. Finally, we also provide a discussion how to use the current results to address the multi sensor/link failure scenario.

This paper is organized as follows. The problem statements are presented in Section

2. Section 3 introduces concepts of structural systems theory and some known results. In Section 4, we present the main results (the proofs are relegated to Appendix). In Section 5, an illustrative example is provided. Finally, conclusions and further research directions are presented in Section 6.

## 2   PROBLEM FORMULATION

Let $\bar{A} \in \{0, \star\}^{n \times n}$ and $\bar{\mathbb{I}}_n(\mathcal{J}) \in \{0, \star\}^{|\mathcal{J}| \times n}$ be the *structural pattern* (i.e., location of fixed zeros or independent parameters) of the system dynamics in (1)-(2). With some common abuse of terminology, we refer to those entries representing independent parameters as "non-zero", while bearing in mind that its realization can be any real number, including zero. In addition, denote by $[\bar{A}]_{i,j}$ the entry in the $i$-th row and $j$-th column of matrix $\bar{A}$, where $[\bar{A}]_{i,j} = 0$ if all possible numerical realizations of the physical system have a fixed zero, and $[\bar{A}]_{i,j} = \star$ otherwise. Let (1)-(2) be described by the pair $(A, \bar{\mathbb{I}}_n(\mathcal{J}))$. Then, a pair $(\bar{A}, \bar{\mathbb{I}}_n(\mathcal{J}))$ is said to be structurally observable if there exists an observable pair $(A', \bar{\mathbb{I}}_n(\mathcal{J}))$ with the same structure (i.e., same locations of zeroes and non-zeroes) as $(\bar{A}, \bar{\mathbb{I}}_n(\mathcal{J}))$. In fact, if a pair $(\bar{A}, \bar{\mathbb{I}}_n(\mathcal{J}))$ is structurally observable, then almost all pairs with the same structure are observable [27].

In this paper, we are interested in the robust output design problems $\mathbf{P_s}$ and $\mathbf{P_l}$, which in the sequel will be formally introduced as $\mathcal{P}_s$ and $\mathcal{P}_l$, respectively.

$\mathcal{P}_s$ Robustness with respect to sensor failure: Given $\bar{A}$, determine a minimal subset of dedicated sensors $\mathcal{J}^* \subset \{1, \ldots, n\}$ such that $(\bar{A}, \bar{\mathbb{I}}_n(\mathcal{J}^*))$ is structurally observable with respect to any single dedicated sensor failure.

In the paper, we will mainly focus on the case where the links are considered as directed, and extend the results to the case where we have undirected link failures in Section 4.2.2. Now, a (directed) link in the dynamic matrix is associated with a non-zero entry $[\bar{A}]_{i,j} = \star$, with graphical representation formalized in Section 4.2. Therefore, a directed link failure corresponds to changing a non-zero entry $[\bar{A}]_{i,j}$ in the structural dynamic matrix $\bar{A}$ by a zero entry, while the rest of elements in $\bar{A}$ remains the same. We denote the new structural dynamic matrix after the link failure by $\bar{A}_{-(i,j)}$.

Therefore the sensor placement ensuring robustness with respect to a link failure can be stated as follows.

$\mathcal{P}_l$ Robustness with respect to link failure: Given $\bar{A}$, determine a minimal subset of dedicated sensors $\mathcal{J}^* \subset \{1, \ldots, n\}$, such that the pair $(\bar{A}_{-(i,j)}, \bar{\mathbb{I}}_n(\mathcal{J}^*))$ is structurally observable when the failure of an arbitrary direct link $[\bar{A}]_{i,j}$ occurs.

Note that when $\mathcal{J}^*$ is obtained as described above, the pair $(\bar{A}, \bar{\mathbb{I}}_n(\mathcal{J}^*))$, i.e., the nominal system without a link failure, is structurally observable.

In addition, we explore the scenario where multiple sensor/link failures can occur, and discuss how the proposed solution to the above problems can be extended to cope with the former scenario.

## 3 PRELIMINARIES AND TERMINOLOGY

A desirable aspect of structural systems is that their properties, such as structural observability, can be analyzed by means of graph theoretical tools.

Let a digraph (i.e., directed graph) be expressed as $\mathcal{D} = (\mathcal{V}, \mathcal{E})$, in which $\mathcal{V}$ denotes a set of *vertices* and $\mathcal{E}$ represents a set of *edges*, such that, an edge $(v_j, v_i)$ is directed from vertex $v_j$ to vertex $v_i$. In addition, we denote by $|\mathcal{V}|$ the number of elements in the set $\mathcal{V}$.

A *path* is a digraph $\mathcal{P} = (\mathcal{V}, \mathcal{E})$ with $\mathcal{V} = \{v_1, \cdots, v_n\}$ and $\mathcal{E} = \{(v_1, v_2), (v_2, v_3), \cdots, (v_{n-1}, v_n)\}$. Further, $v_1$ and $v_n$ are referred to as *root* and *tip* of the path, respectively, and an *isolated vertex* (with some abuse of notation) is also considered as a (degenerated) simple path. A cycle can be defined as a path with an additional edge from its tip to its root.

Given a digraph $\mathcal{D} = (\mathcal{V}, \mathcal{E})$, $\mathcal{D}_S = (\mathcal{V}_S, \mathcal{E}_S)$ is a *subgraph* of $\mathcal{D}$ if $\mathcal{V}_S \subset \mathcal{V}$ and $\mathcal{E}_S \subset \mathcal{E}$. A digraph $\mathcal{D}$ is said to be strongly connected if there exists a path between any pair of vertices [28]. A *strongly connected component* (SCC) is a maximal subgraph $\mathcal{D}_S = (\mathcal{V}_S, \mathcal{E}_S)$ of $\mathcal{D}$, i.e., it contains the largest collection of vertices and edges between these, such that for every $v, w \in \mathcal{V}_S$ there exists a path from $v$ to $w$. Any digraph $\mathcal{D} = (\mathcal{V}, \mathcal{E})$ can be uniquely decomposed into disjoint SCCs. Visualizing each SCC as a supernode, one can generate a *directed acyclic graph* (DAG), i.e., a digraph with no cycles, in which each supernode corresponds to an SCC and an edge exists between two SCCs if and only if there exists an edge between vertices belonging to the corresponding SCCs. The unique DAG associated with $\mathcal{D} = (\mathcal{V}, \mathcal{E})$ can be efficiently generated in $\mathcal{O}(|\mathcal{V}| + |\mathcal{E}|)$ [28], where $|\mathcal{V}|$ and $|\mathcal{E}|$ denote the number of vertices and edges in $\mathcal{D}$. We can further characterize the SCCs in the DAG representation of $\mathcal{D}$ as follows: an SCC is a *sink-SCC* if it has no outgoing edges from its vertices to those in another SCC.

A digraph $\mathcal{D} = (\mathcal{V}, \mathcal{E})$ can be spanned by a disjoint union of paths and cycles. Therefore, we introduce the following decomposition of a digraph.

**Definition 1.** *A path and cycle (P&C) decomposition of a digraph $\mathcal{D} = (\mathcal{V}, \mathcal{E})$ consists of a disjoint union of paths and cycles that span $\mathcal{D}$. In addition, a* minimum P&C decomposition *is a P&C decomposition with the minimum number of paths among possible P&C decompositions that span $\mathcal{D}$.* □

The minimum P&C decomposition is useful in explaining some preliminaries and the main results of the paper. Note that the minimum P&C decomposition of a digraph is generally not unique, while the total number of paths in all minimum P&C decompositions is unique. Further, a minimum P&C decomposition of the state digraph can be obtained efficiently in $\mathcal{O}(\sqrt{|\mathcal{X}|}|\mathcal{E}|)$ [6].

Now, consider the structure associated to the dynamical system plant (1) and (2), given by $\bar{A}$ and $\bar{\mathbb{I}}_n(\mathcal{J})$ as discussed in Section 2. We can associate with each state variable a vertex in a digraph, to which we refer to as a *state vertex* and its collection is given by $\mathcal{X} = \{x_1, \cdots, x_n\}$. In addition, the edges between state vertices are related with non-zero entries of the matrix $\bar{A}$, that we refer to as *state edges*, and we denoted by $\mathcal{E}_{\bar{A}} = \{(x_i, x_j) : [\bar{A}]_{ji} \neq 0\}$. Subsequently, we define the *state digraph* as $\mathcal{D}(\bar{A}) = (\mathcal{X}, \mathcal{E}_{\bar{A}})$, used hereafter to study structural observability.

The proposed solutions in Section 4 consist of recasting $\mathcal{P}_s$ and $\mathcal{P}_l$ as *weighted set covering problems* [29], that may be described as follows: consider a set $\mathcal{U}$, called the universal set, and $k$ sets $\mathcal{C}_i \subset \mathcal{U}$, with $i \in \mathcal{I} = \{1, \cdots, k\}$. Each $\mathcal{C}_i$ is associated with a cost $c_i \in \mathbb{R}_0^+$. The goal is to find a collection $\{\mathcal{C}_i\}_{i \in \mathcal{J}^*}$, with $\mathcal{J}^* \subset \mathcal{I}$, such that

$$\mathcal{J}^* = \arg \min_{\mathcal{J} \subset \mathcal{I}} \quad \sum_{k \in \mathcal{J}} c_k$$

$$\text{s.t.} \quad \mathcal{U} \subset \bigcup_{j \in \mathcal{J}} \mathcal{C}_j. \tag{3}$$

If the costs are positive and uniform, then we obtain the well known *set covering problem*. Although, the (weighted) set covering problem is NP-hard, some approximation algorithms with polynomial complexity are known. Further, these algorithms are known to ensure bounded optimality gap, hence, providing optimality guarantees.

The collection of state variables $\{x_j\}_{j \in \mathcal{J}}$ measured by dedicated sensors indexed by $\mathcal{J}$ is referred to as a *solution* in this paper. Next, we introduce the concept of *feasible solution* as follows.

**Definition 2** ([6]). *Let* $\mathcal{D}(\bar{A}) = (\mathcal{X} \equiv \{x_1, \ldots, x_n\}, \mathcal{E}_{\bar{A}})$ *be the state digraph. A solution* $\mathcal{F} = \cup_{j \in \mathcal{J}} \{x_j\}$ *is a* feasible *solution if* $(\bar{A}, \bar{\mathbb{I}}_n(\mathcal{J}))$ *is structurally observable, where the letter* $\mathcal{F}$ *refers to feasible.* □

Now, we revisit a characterization of the feasible solutions.

**Theorem 1** ([6]). *Let* $\mathcal{D}(\bar{A}) = (\mathcal{X}, \mathcal{E}_{\bar{A}})$ *denote the state digraph. A set* $\mathcal{F} \subset \mathcal{X}$ *is a feasible solution if and only if there exist two subsets* $\mathcal{T}$ *and* $\mathcal{S}^\circ$ *such that* $\mathcal{F} \supset \mathcal{T} \cup \mathcal{S}^\circ$, *where* $\mathcal{T}$ *contains the tips of the paths in a minimum P&C decomposition, and* $\mathcal{S}^\circ$ *denotes a subset containing at least one variable from each sink-SCC that does not contain variables from* $\mathcal{T}$; *in particular, each of the sink-SCCs that does not contain variables from* $\mathcal{T}$ *is spanned by cycles.* □

In addition, a feasible solution with the minimal number of state variables is said to be a *minimal feasible solution*. Next, we present a characterization of the minimal feasible solutions.

**Theorem 2** ([6]). *Let* $\mathcal{D}(\bar{A}) = (\mathcal{X}, \mathcal{E}_{\bar{A}})$ *denote the state digraph. A set* $\mathcal{F} \subset \mathcal{X}$ *is a minimal feasible solution if and only if it satisfies the condition of Theorem 1, and* $\mathcal{T}$ *has the maximum number of tips of the paths in as many distinct sink-SCCs among all possible minimum P&C decompositions.* □

Recall that the minimum P&C decomposition may not be unique, hence, although they have the same number of paths, the tips of the paths associated with a minimum P&C decomposition change. In particular, they may belong to distinct SCCs of the state digraph; in particular, the tips of the paths may or may not be in a sink-SCC. Note that a sink-SCC that contains at least one tip of paths will not contribute to the size of any feasible solution. Consequently, if the tip of a path belongs to a sink-SCC, we say that this state vertex plays a *double role* since it simultaneously satisfies the two conditions in Theorem 1 and 2. In order to determine a minimal feasible solution, i.e., the feasible solution with the lowest cardinal among these, we just need to maximize the number of state vertices with double role, i.e., that are both tips of the paths and lie in different sink-SCCs.

For example, consider a $10$-state system (randomly generated), whose state digraph is depicted in Fig. 5-(a). The state digraph consists of four SCCs, as depicted in the four rectangles in Fig. 5-(b). Note that the only sink-SCC is $\{x_1, x_2, x_4, x_5, x_7, x_8, x_{10}\}$, and is enclosed by the top rectangle. In addition, a minimum P&C decomposition of the state digraph is also depicted in the black continuous arrows and vertices in Fig. 5-(b), which consists of two paths and one cycle. the tips of the paths associated with the given minimum P&C decomposition of the state digraph is $\mathcal{T} = \{x_8, x_{10}\}$. Hence, from Theorem 2, a minimal feasible solution $\mathcal{F}$ contains the tips of the paths $\mathcal{T}$ associated with the minimum P&C decomposition of the state digraph. In addition, the only sink-SCC contains elements $x_8$ and $x_{10}$ in $\mathcal{T}$; hence, $\mathcal{F} = \mathcal{T} = \{x_8, x_{10}\}$ is a minimal feasible solution, demonstrated by grey squares in Fig. 5-(b).

In [6], it was further shown that a minimal feasible solution can be determined using polynomial complexity algorithms in the size of the state space.

## 4 MAIN RESULTS

In this section, we present the main results of the paper, i.e., we address the problems formulated in $\mathcal{P}_s$ and $\mathcal{P}_l$. In Section 4.1, we address $\mathcal{P}_s$, followed by $\mathcal{P}_l$ addressed in Section 4.2. Briefly, $\mathcal{P}_s$ and $\mathcal{P}_l$ are first re-casted in terms of feasible solutions (see Definition 2) in $\mathcal{P}_s'$ and $\mathcal{P}_l'$, respectively. Subsequently, they are shown to be NP-hard in Theorem 3 and Theorem 6, respectively. Due to the computational hardness of these problems, an intuitive two-step approach to $\mathcal{P}'_s$ and $\mathcal{P}'_l$ is provided, which ultimately reduces to a weighted set covering problem. Feasibility properties of a solution to $\mathcal{P}'_s$ using the proposed approach is provided in Theorem 4, and its polynomial-time approximation in Theorem 5. Similarly, feasibility properties of a solution to $\mathcal{P}'_l$ is stated in Theorem 7, and its polynomial-time approximation in Theorem 8. Finally, in Section 4.3 we explain how the proposed solutions can be extended to deal with the multiple sensor/link failures scenario.

## 4.1  Robustness with respect to sensor failure

First, we recast $\mathcal{P}_s$ in terms of feasible solutions as follows.

$\mathcal{P}'_s$: Given the state digraph $\mathcal{D}(\bar{A}) = (\mathcal{X}, \mathcal{E}_{\bar{A}})$, select a minimum subset of state variables $\mathcal{F}^s \subset \mathcal{X}$, such that for all $x \in \mathcal{F}^s$, $\mathcal{F}^s \setminus \{x\}$ is a feasible solution. ∘

A feasible solution $\mathcal{F}^s$ with the properties described in $\mathcal{P}'_s$ is referred to as a *sensor-robust feasible solution (s-robust feasible solution)*. In addition, to ease the analysis, when referring to $\mathcal{F}^s \setminus \{x\}$ we say that $\{x\}$ is discarded from $\mathcal{F}^s$, and the remaining variables $\mathcal{F}^s \setminus \{x\}$ are said to be non-discarded.

Now, we show the hardness of $\mathcal{P}'_s$ (or, equivalently, $\mathcal{P}_s$).

**Theorem 3.** *$\mathcal{P}'_s$ is NP-hard.* □

Given the above hardness, we provide an approximation approach to $\mathcal{P}'_s$. To determine an $s$-robust feasible solution, we recall Theorem 2. More precisely, consider a state digraph $\mathcal{D}(\bar{A}) = (\mathcal{X}, \mathcal{E}_{\bar{A}})$ and one associated minimal feasible solution $\mathcal{F}(\Xi) = \{x_{\tau_1}, \cdots, x_{\tau_p}\}$, where $\tau_i \in \{1, \cdots, n\}$ (for $i \in \{1, \cdots, p\}$) is the index of each state variable, and $\Xi$ explicitly states the minimum P&C decomposition considered (see Theorem 1 and 2). More specifically, $\mathcal{F}(\Xi)$ contains two subsets $\mathcal{T}(\Xi)$ and $\mathcal{S}^{\circ}(\Xi)$, where (without loss of generality) $\mathcal{T}(\Xi) = \{x_{\tau_1}, \cdots, x_{\tau_m}\}$ corresponds to the tips of the paths in $\Xi$, and $\mathcal{S}^{\circ}(\Xi) = \{x_{\tau_{m+1}}, \cdots, x_{\tau_p}\}$ is the set of state vertices chosen seriatim from each sink-SCC without tips from $\mathcal{T}(\Xi)$.

Now, observe that a discarded state variable can belong to either $\mathcal{T}(\Xi)$ or $\mathcal{S}^\circ(\Xi)$. Thus, we need to consider new state variables that together with the non-discarded state variables from the minimal feasible solution originates an $s$-robust feasible solution. To sharpen the intuition on the two-step algorithm presented hereafter, in Example 1 and 2 we provide examples of discarded state variables that belong to $\mathcal{S}^\circ(\Xi)$ and $\mathcal{T}(\Xi)$, and explain how to construct an $s$-robust feasible solution using a minimal feasible solution.

**Example 1.** *Let $\mathcal{D}(\bar{A}) = (\mathcal{X}, \mathcal{E}_{\bar{A}})$ consist of a single SCC and have a minimum P&C decomposition consisting of only cycles. From Theorem 2, we have that $\mathcal{F} = \{x\}$ with any $x \in \mathcal{X}$, is a minimal feasible solution. Subsequently, $\mathcal{F}^s = \mathcal{F} \cup \{x'\}$, with $x' \in \mathcal{X} \setminus \{x\}$, is a (minimal) $s$-robust feasible solution.* $\qquad\square$

We refer to the additional state variables, that together with those in $\mathcal{F}$ lead to an $s$-robust feasible solution, as *alternatives* to $\mathcal{F}$. This leads us to compute alternatives for an arbitrary minimal feasible solution, which correspond to the state variables $\mathcal{F}'$, such that after any $x \in \mathcal{F}$ is discarded, $(\mathcal{F} \setminus \{x\}) \cup \mathcal{F}'$ still forms an $s$-robust feasible solution. Further, recall that a minimal feasible solution $\mathcal{F}(\Xi)$ contains sets $\mathcal{T}(\Xi)$ and $\mathcal{S}^\circ(\Xi)$; subsequently, let us start by introducing the alternatives for elements in $\mathcal{S}^\circ$, as considered in Example 1, and define the set of *sink-alternatives* $\Delta_{\tau_i}^{\mathcal{S}}$ to a particular state vertex $x_{\tau_i} \in \mathcal{S}^\circ(\Xi)$.

**Definition 3** (Sink-alternatives). *Consider a state digraph $\mathcal{D}(\bar{A}) = (\mathcal{X}, \mathcal{E}_{\bar{A}})$ and a minimal feasible solution $\mathcal{F}(\Xi) = \mathcal{T}(\Xi) \cup \mathcal{S}^\circ(\Xi)$, where $\mathcal{F}(\Xi) = \{x_{\tau_1}, \cdots, x_{\tau_p}\}$, and (without loss of generality) $\mathcal{T}(\Xi) = \{x_{\tau_1}, \cdots, x_{\tau_m}\}$, and $\mathcal{S}^\circ(\Xi) = \{x_{\tau_{(m+1)}}, \cdots, x_{\tau_p}\}$. Let $n_{(i-m)}$ denote the number of state variables in the same sink-SCC as $x_{\tau_i}$ (not including $x_{\tau_i}$), where $i > m$. Then,*

$$\Delta_{\tau_i}^{\mathcal{S}} = \Big\{ x \in (\mathcal{X} \setminus \{x_{\tau_i}\}) : x \text{ belongs to the same sink-SCC as } x_{\tau_i} \Big\}, \qquad (4)$$

*is the set of sink-alternatives to state vertex $x_{\tau_i} \in \mathcal{S}^\circ(\Xi)$. Further, for notational convenience, we denote the state variables in $\Delta_{\tau_i}^{\mathcal{S}}$ as $\Delta_{\tau_i}^{\mathcal{S}} = \{\delta_{\tau_i,1}^{\mathcal{S}}, \cdots, \delta_{\tau_i,n_{(i-m)}}^{\mathcal{S}}\}$.* $\qquad\square$

Similarly, we can now introduce the notion of alternatives for elements in $\mathcal{T}$, which we briefly refer to as *tip-alternatives*, motivated by the following example.

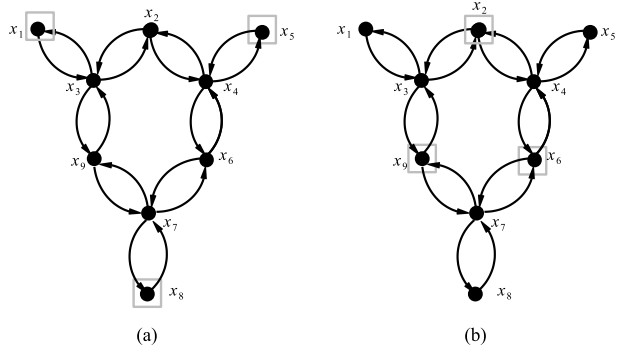

Fig. 1. An example of a state digraph $\mathcal{D}(\bar{A})$ with 9 state variables. In (a) and (b) we depict two minimal feasible solutions $\mathcal{F}_1 = \{x_1, x_5, x_8\}$ and $\mathcal{F}_2 = \{x_2, x_6, x_9\}$ enclosed by squares, respectively. In particular, the two minimal feasible solutions lead to $s$-robust feasible solutions with different sizes based on our approach.

**Example 2.** *Fig. 1-(a) depicts a $9$-state system and a minimal feasible solution $\mathcal{F}(\Xi) = \mathcal{T}(\Xi) \cup \mathcal{S}^\circ(\Xi)$, where $\mathcal{T}(\Xi) = \{x_1, x_5, x_8\}$, and $\mathcal{S}^\circ(\Xi) = \emptyset$, because the state digraph consists of one sink-SCC that contains tips of the paths. Discarding either $x_1, x_5$ or $x_8$ from $\mathcal{F}(\Xi)$ renders a non-feasible solution, i.e., structural observability no longer holds. As discussed in Section 3, the minimum P&C decomposition of the state digraph is not unique. In fact, $\mathcal{T}'(\Xi') = (\mathcal{T}(\Xi) \setminus \{x\}) \cup \{y\}$, with $x \in \mathcal{T}(\Xi)$ and $y \in \{x_2, x_6, x_9\}$, is also the tips of the paths associated with another minimum P&C decomposition $\Xi'$. Hence, $\mathcal{F}(\Xi) \cup \{y\}$, with $y \in \{x_2, x_6, x_9\}$, is an $s$-robust feasible solution. Furthermore, this $s$-robust feasible solution is also minimal.* □

Subsequently, we define the set of *tip-alternatives* $\Delta_{\tau_i}^{\mathcal{T}}$ to a state variable $x_{\tau_i} \in \mathcal{T}$ as follows.

**Definition 4** (Tip-alternatives)**.** *Consider a state digraph $\mathcal{D}(\bar{A}) = (\mathcal{X}, \mathcal{E}_{\bar{A}})$ and a minimal feasible solution $\mathcal{F}(\Xi) = \mathcal{T}(\Xi) \cup \mathcal{S}^\circ(\Xi)$, where $\mathcal{F}(\Xi) = \{x_{\tau_1}, \cdots, x_{\tau_p}\}$, and (without loss of generality) $\mathcal{T}(\Xi) = \{x_{\tau_1}, \cdots, x_{\tau_m}\}$, and $\mathcal{S}^\circ(\Xi) = \{x_{\tau_{(m+1)}}, \cdots, x_{\tau_p}\}$. Then, for all $i \leq m$,*

$$\Delta_{\tau_i}^{\mathcal{T}} = \Big\{ x \in (\mathcal{X} \setminus \{x_{\tau_i}\}) : (\mathcal{T}(\Xi) \setminus \{x_{\tau_i}\}) \cup \{x\} \text{ is a set of tips of}$$
$$\text{paths associated with a minimum P\&C decomposition of } \mathcal{D}(\bar{A}) \Big\} \quad (5)$$

*is the set of tip-alternatives to a tip of a path $x_{\tau_i} \in \mathcal{T}(\Xi)$. Further, for notational convenience,*

we denote the state variables in $\Delta_{\tau_i}^{\mathcal{T}}$ as $\Delta_{\tau_i}^{\mathcal{T}} = \{\delta_{\tau_i,1}^{\mathcal{T}}, \cdots, \delta_{\tau_i,\mu_i}^{\mathcal{T}}\}$, where $\mu_i$ is the number of tip-alternatives to $x_{\tau_i}$. □

In other words, keeping all state variables in the tips of the paths $\mathcal{T}$ fixed except for $x_{\tau_i}$ that is replaced by a state variable in $\Delta_{\tau_i}^{\mathcal{T}}$, we obtain the tips of the paths associated with another minimum P&C decomposition $\Xi'$. For example, in Example 2, we have $\Delta_1^{\mathcal{T}} = \Delta_5^{\mathcal{T}} = \Delta_8^{\mathcal{T}} = \{x_2, x_6, x_9\}$, i.e., $x_2$, $x_6$, and $x_9$ are tip-alternatives to either $x_1$, $x_5$, or $x_8$.

Now, notice that given a minimal feasible solution $\mathcal{F} = \mathcal{T} \cup \mathcal{S}^\circ$, replacing a state variable $x_{\tau_i} \in \mathcal{S}^\circ$ by a state variable in its sink-alternatives $\Delta_{\tau_i}^{\mathcal{S}}$, i.e., $\mathcal{F}' = \mathcal{T} \cup (\mathcal{S}^\circ \setminus \{x_{\tau_i}\} \cup \{x\})$ with $x \in \Delta_{\tau_i}^{\mathcal{S}}$, $\mathcal{F}'$ is still a minimal feasible solution. Alternatively, when tip-alternatives are considered two situations can occur. For instance, in Example 1 feasibility is retained mainly due to the existence of a single SCC. However, replacing a state variable $x_{\tau_i} \in \mathcal{T}$ by one of its tip-alternatives may not yield a feasible solution in general. More precisely, for a tip $x_{\tau_i} \in \mathcal{T}$, whether $\mathcal{F}_x^{-x_{\tau_i}} = \{\mathcal{T} \setminus \{x_{\tau_i}\} \cup \{x\}\} \cup \mathcal{S}^\circ$ with $x \in \Delta_{\tau_i}^{\mathcal{T}}$, which we refer to as an *interchanged solution*, is a feasible solution or not, depends on which SCCs of $\mathcal{D}(\bar{A})$ the variables $x_{\tau_i}$ and $x$ belong to. In particular, if in an interchanged solution $\mathcal{F}_x^{-x_{\tau_i}}$, $x_{\tau_i}$ is the only tip of paths in a sink-SCC, by replacing $x_{\tau_i}$ by one of its tip-alternatives that does not belong to the same SCC, this SCC becomes an SCC without a tip of the paths in a minimum P&C decomposition; thus, a new state variable in this SCC needs to be considered to satisfy the second requirement in Theorem 2. Hence, we need to replace $x_{\tau_i}$ by two state variables: a tip-alternative $x$, and a state variable $x'$ (different from $x_{\tau_i}$) in the same sink-SCC as $x_{\tau_i}$.

In summary, we can generalize our previous discussion by considering collection of sets that account for the sink- and tip-alternatives for each variable $x_{\tau_i} \in \mathcal{F}$, which we refer to as *back-ups* and which we denote by $\Omega_{\mathcal{F}}^i$. In particular, the sets in $\Omega_{\mathcal{F}}^i$ can contain one state variable, when sink-alternatives are considered, which we denote by $\{\delta_{\tau_i,j}^{\mathcal{S}}\}$, or tip-alternatives that do not increase the number of sink-SCCs without tips of the paths, which we denote by $\omega_{\mathcal{F}}^{ijk}$. However, when tip-alternatives increase the number of sink-SCCs without tip of paths, then $\omega_{\mathcal{F}}^{ijk}$ contains two state variables: one is a tip-alternative and another is a state variable in the same sink-SCC as the discarded state variable. Subsequently, we have the following definition.

**Definition 5** (Set of Back-ups). *Consider a state digraph $\mathcal{D}(\bar{A}) = (\mathcal{X}, \mathcal{E}_{\bar{A}})$ and a minimal feasible solution $\mathcal{F} = \mathcal{T} \cup \mathcal{S}^\circ$, where $\mathcal{F} = \{x_{\tau_1}, \cdots, x_{\tau_p}\}$, and (without loss of generality) $\mathcal{T} = \{x_{\tau_1}, \cdots, x_{\tau_m}\}$, and $\mathcal{S}^\circ = \{x_{\tau_{(m+1)}}, \cdots, x_{\tau_p}\}$. Then, the set of back-ups to $x_{\tau_i}$, referred to as $\Omega_{\mathcal{F}}^i$, is given by*

1) *for $i = 1, \ldots, m$,*

$$\Omega_{\mathcal{F}}^i = \{\Omega_{\mathcal{F}}^{ij}\}_{j \in \{1, \cdots, \mu_i\}} = \{\omega_{\mathcal{F}}^{ijk}\}_{j \in \{1, \cdots, \mu_i\}, k \in \{1, \cdots, n_{i,j}\}}, \tag{6}$$

*where $\mu_i$ is the total number of tip-alternatives of $x_{\tau_i}$ (see Definition 4), $n_{i,j}$ equals to the total number of state variables that belong to the same sink-SCC as $x_{\tau_i}$ but not $x_{\tau_i}$, $\omega_{\mathcal{S}}^{ijk} = \{\delta_{\tau_{i,j}}^{\mathcal{T}}\} \cup \mathcal{A}_{i,j,k}$, with $\mathcal{A}_{i,j,k} \equiv \emptyset$ if the interchanged solution $\mathcal{F}_{\delta_{\tau_{i,j}}^{\mathcal{T}}}^{-x_{\tau_i}}$ is a feasible solution, and $\mathcal{A}_{i,j,k} = \{x_{i_k}\}$ otherwise, where $x_{i_k} \neq x_{\tau_i}$ is the k-th state variable that belongs to the same SCC as $x_{\tau_i}$; and*

2) *for $i = m+1, \ldots, p$,*

$$\Omega_{\mathcal{F}}^i = \{\Omega_{\mathcal{F}}^{ij}\}_{j \in \{1, \cdots, n_{(i-m)}\}} = \{\{\delta_{\tau_{i,j}}^{\mathcal{S}}\}\}_{j = \{1, \cdots, n_{(i-m)}\}}, \tag{7}$$

*where $n_{(i-m)}$ is the total number of sink-alternatives of $x_{\tau_i}$ (see Definition 3).* $\square$

Therefore, the goal is to find a minimal $\mathcal{F}'$ with above properties, i.e., there is no other $\mathcal{F}''$ with $|\mathcal{F}''| < |\mathcal{F}'|$ such that $\mathcal{F} \cup \mathcal{F}''$ is an $s$-robust feasible solution. Notice that a common back-up may exist for different elements in $\mathcal{F}$. For example, in Example 1, $\{x_2\}$ is a common back-up for $x_1, x_5$ and $x_8$. As a result, minimizing the size of $\mathcal{F}'$ is the same as maximizing the shared state variables in back-up sets for different elements in $\mathcal{F}$. Subsequently, we can find the minimal $\mathcal{F}'$ by considering a set covering where the sets account for the back-ups for any element in $\mathcal{F}$. To this goal, consider the construction of back-ups in (6)-(7). Since $\Omega_{\mathcal{F}}^i$ in (6)-(7) consist of sets with one or two state variables, to avoid cumbersome notation, let $\mathcal{Z}$ be a collection of all $N = n + \binom{n}{2} = \frac{n(n+1)}{2}$ possible sets denoted by $\mathcal{Z} = \{\mathcal{Z}_1, \cdots, \mathcal{Z}_n, \mathcal{Z}_{n+1}, \cdots, \mathcal{Z}_N\}$, where $\mathcal{Z}_i = \{x_i\}$ (for $i = 1, \cdots, n$), $\mathcal{Z}_i =$ $\{x_\alpha, x_\beta\}\ \forall \alpha, \beta \in \{1, \cdots, n\}$ $(\alpha \leq \beta)$, and $i = \pi(\alpha-1, \beta-1) = \frac{1}{2}(\alpha+\beta-2)(\alpha+\beta-1)+\beta-1+n$

denotes the Cantor pairing function shifted by $n$. In summary,

$$\mathcal{Z} = \{\{x_1\}, \cdots, \{x_n\}, \{x_1, x_2\}, \cdots, \{x_{n-1}, x_n\}\}. \tag{8}$$

Consequently, from (6)-(7), we can see that any element in $\Omega_{\mathcal{F}}^i$ consists of one element in $\mathcal{Z}$, which implies that $\Omega_{\mathcal{F}}^i = \{\mathcal{Z}_j\}_{j \in \mathcal{Q}}$, with $\mathcal{Q} \subset \{1, \ldots, N\}$. With this definition, the problem reduces to identifying for each element $\mathcal{Z}_i$, the subset of state variables of $\mathcal{F}$ for which $\mathcal{Z}_i$ is a shared back-up. This can be achieved by considering the definition of the sets $\mathcal{V}_j^s$ that contain the indices of the state variables in a minimal feasible solution $\mathcal{F}$, or equivalently $\Omega_{\mathcal{F}}^i$, that an element $\mathcal{Z}_j$ is a back-up to. Formally, this can be defined as follows: let $\mathcal{Z} = \{\mathcal{Z}_1, \cdots, \mathcal{Z}_N\}$ and $\mathcal{I} = \{1, \cdots, p\}$, then

1)  for $j = 1, \cdots, n$,

$$\mathcal{V}_j^s = \{i \in \mathcal{I} : \mathcal{Z}_j \in \Omega_{\mathcal{F}}^i\}, \tag{9}$$

2)  and for $j = n+1, \cdots, N$,

$$\mathcal{V}_j^s = \{i \in \mathcal{I} : \mathcal{Z}_j \in \Omega_{\mathcal{F}}^i\} \cup \left\{ \bigcup_{x_k \in \mathcal{Z}_j} \mathcal{V}_k^s \right\}. \tag{10}$$

Therefore, if $\mathcal{Z}_j$ consists of two state variables, then $\mathcal{Z}_j$ is not only a back-up of those $x_{\tau_i}$'s satisfying $\mathcal{Z}_j \in \Omega_{\mathcal{F}}^i$, but also the back-ups of those $x_{\tau_i}$'s to which the back-ups contain an individual state variable belonging to $\mathcal{Z}_j$. To sum up, by computing $\mathcal{V}_i^s$, we identify the subset of state variables of $\mathcal{F}$ for which $\mathcal{Z}_i \in \mathcal{Z}$ is a shared back-up.

Thus, we can determine an $s$-robust feasible solution as follows:

**Lemma 1.** *Consider a state digraph $\mathcal{D}(\bar{A}) = (\mathcal{X}, \mathcal{E}_{\bar{A}})$ and one of its minimal feasible solutions $\mathcal{F} = \{x_{\tau_1}, x_{\tau_2}, \cdots, x_{\tau_p}\}$. Consider the sets in (9)-(10), and let $\mathcal{I} = \{1, \cdots, p\}$. If there exists $\mathcal{J} \subset \{1, \cdots, N\}$ such that $\mathcal{I} \subset \bigcup_{j \in \mathcal{J}} \mathcal{V}_j^s$, i.e., the family $\{\mathcal{V}_j^s\}_{j \in \mathcal{J}}$ covers $\mathcal{I}$, then*

$$\mathcal{F}^* = \bigcup_{j \in \mathcal{J}} \mathcal{Z}_j \ \cup \ \mathcal{F},$$

*is an $s$-robust feasible solution.* $\qquad \square$

 

Notice that Lemma 1 only determines the minimum number of back-ups based on a particular seed feasible solution $\mathcal{F}$. Furthermore, if some $\mathcal{Z}_j$ with two state variables is considered in the covering, there might exist a set $\mathcal{Z}_{j'}$ that consists of a single state variable that could have been used instead to cover $\mathcal{I}$. Thus, the minimum number of back-ups does not ensure the minimum number of alternatives to be considered (in addition to a minimal feasible solution), which motivates us to introduce a weighed version of the problem to account for this scenario. More precisely, given the sets in (9)-(10), we assign to each $\mathcal{V}_i^s$ a cost $c_i$, given by

$$c_i = |\mathcal{Z}_i|, \text{ for } i = 1, \cdots, N. \tag{11}$$

In other words, the weight of $\mathcal{V}_i^s$ for the set covering problem is equal to the number of state vertices in $\mathcal{Z}_i$. Hence, we obtain the following result:

**Theorem 4.** *Consider a state digraph $\mathcal{D}(\bar{A}) = (\mathcal{X}, \mathcal{E}_{\bar{A}})$ and one of its minimal feasible solutions $\mathcal{F} = \{x_{\tau_1}, x_{\tau_2}, \cdots, x_{\tau_p}\}$. Consider the sets (9)-(11), and let $\mathcal{I} = \{1, \cdots, p\}$. If there exists $\mathcal{J}^*$ satisfying:*

$$\begin{aligned} \mathcal{J}^* = \underset{\mathcal{J} \subset \{1, \cdots, N\}}{\arg \min} \quad & \sum_{k \in \mathcal{J}} c_k, \\ \text{subject to} \quad & \mathcal{I} \subset \bigcup_{j \in \mathcal{J}^*} \mathcal{V}_j^s, \end{aligned} \tag{12}$$

*i.e., the family $\{\mathcal{V}_j^s\}_{j \in \mathcal{J}^*}$ covers $\mathcal{I}$, then*

$$\mathcal{F}^* = \bigcup_{j \in \mathcal{J}^*} \mathcal{Z}_j \ \cup \ \mathcal{F},$$

*is an $s$-robust feasible solution. In addition, there is no other $s$-robust feasible solution $\mathcal{F}'$ with $\mathcal{F} \subset \mathcal{F}'$ satisfying Lemma 1 such that $|\mathcal{F}'| < |\mathcal{F}^*|$.* □

Theorem 4 requires to solve an NP-hard problem – the weighted set covering problem [28]. To circumvent that problem, the next result shows how to obtain an approximated solution, which can be achieved using polynomial complexity algorithms [30], [31] with optimality guarantees.

**Theorem 5.** *Consider a state digraph $\mathcal{D}(\bar{A}) = (\mathcal{X}, \mathcal{E}_{\bar{A}})$, an arbitrary minimal feasible solution $\mathcal{F}$, the sets (9)-(11), and the weighted set covering problem as described in Theorem 4. Then,* *the corresponding weighted set covering can be constructed with complexity $\mathcal{O}(|\mathcal{X}|^5)$. Further, a feasible but approximate solution $\mathcal{J}'$ of the corresponding minimum weighted set covering may be constructed using an $\mathcal{O}(|\mathcal{X}|^3|)$ complexity algorithm such that*

$$\mathcal{F}' = \bigcup_{j \in \mathcal{J}'} \mathcal{Z}_j \, \cup \, \mathcal{F},$$

*is an $s$-robust feasible solution, and the* performance ratio, *i.e., the ratio between the size of the approximated $s$-robust feasible solution and the size of the minimal $s$-robust feasible solution containing $\mathcal{F}$ is bounded by the harmonic number $H(p) = \sum_{k=1}^{p} \frac{1}{k}$ of $p$, where $p = |\mathcal{F}|$.* □

Finally, we note that the proposed two-step approach only ensures minimality conditioned on a particular initial minimal feasible solution. In other words, the minimality among all possible $s$-robust feasible solutions depends on the minimal feasible solution $\mathcal{F}$ that one starts with. Different minimal feasible solutions may lead to different weighted set covering problems, of which the optimal solutions lead to $s$-robust feasible solutions with different sizes. For instance, consider the example in Fig. 1, which depicts a 9-state system and two of its minimal feasible solutions depicted by grey squares, respectively. On one hand, $\mathcal{F}_1$ leads to alternatives with size one and consequently generates a minimal $s$-robust feasible solution. On the other hand, Fig. 1-(b) depicts another minimal feasible solution $\mathcal{F}_2 = \{x_2, x_6, x_9\}$, and the sets of back-ups for each state vertex in $\mathcal{F}_2$ are $\Omega^1_{\mathcal{F}_2} = \{\{x_1\}, \{x_5\}\}$, $\Omega^2_{\mathcal{F}_2} = \{\{x_5\}, \{x_8\}\}$, and $\Omega^3_{\mathcal{F}_2} = \{\{x_1\}, \{x_8\}\}$, respectively. Hence, by the weighted set covering problem constructed as described in Theorem 4, an $s$-robust feasible solution may be obtained as $\mathcal{F}_2 \cup \{x, y\}$, with $x, y \in \{x_1, x_5, x_8\}$ and $x \neq y$, which is not minimal.

Although our solution does not guarantee a minimal $s$-robust feasible solution, we argue that in practice our set covering based design approach leads to $s$-robust feasible solutions with small number of states variables. For example, consider the system depicted in Fig. 1, for which we can obtain two disjoint minimal feasible solutions $\mathcal{F}_1 = \{x_1, x_5, x_8\}$ and $\mathcal{F}_2 = \{x_2, x_6, x_9\}$. We may then obtain an $s$-robust feasible solution $\mathcal{F}^u_s = \mathcal{F}_1 \cup \mathcal{F}_2$ consisting of 6 state variables, whereas with the proposed approach (Theorem 4), we can obtain an $s$-robust feasible solution consisting of 4 or 5 state variables.

## 4.2 Robustness with respect to link failure

As mentioned in Section 2, in this section, we will mainly focus on the case where the links are considered as directed, and extend the results to the undirected link failure case in Section 4.2.2. The links discussed before Section 4.2.2 are considered as directed links. Similar to Section 4.1, we recast $\mathcal{P}_l$ in terms of feasible solutions. First, we formally define links as follows:

**Definition 6.** (Link) Let $\mathcal{D}(\bar{A}) = (\mathcal{X}, \mathcal{E}_{\bar{A}})$ be the state digraph. Then $\mathcal{L}_{i,j} = \{(x_i, x_j)\} \subset \mathcal{E}_{\bar{A}}$ represents a link from $x_i$ to $x_j$. □

Subsequently, $\mathcal{P}_l$ can be reformulated as follows:

$\mathcal{P}'_l$: Given the state digraph $\mathcal{D}(\bar{A}) = (\mathcal{X}, \mathcal{E}_{\bar{A}})$, select a minimum subset of state variables $\mathcal{F}^l \subset \mathcal{X}$, such that $\mathcal{F}^l$ is a feasible solution for the state digraph $\mathcal{D}(\bar{A}_{-(i,j)}) = (\mathcal{X}, \mathcal{E}_{\bar{A}} - \mathcal{L}_{i,j})$, for all $\mathcal{L}_{i,j} \subset \mathcal{E}_{\bar{A}}$. ○

A feasible solution $\mathcal{F}^l$ with the properties described in $\mathcal{P}'_l$ is referred to as an *link-robust feasible solution (l-robust feasible solution)*.

Further, for convenience, when one link $\mathcal{L}_{i,j}$ fails we say that the state digraph $\mathcal{D}(\bar{A}_{-(i,j)})$ is a *corrupted* state digraph.

Now, we show the hardness of $\mathcal{P}'_l$ (or equivalently $\mathcal{P}_l$).

**Theorem 6.** $\mathcal{P}'_l$ *is NP-hard.* □

### 4.2.1 Directed Link Failures

Similar to Section 4.1, we now provide a procedure to determine an approximate solution to $\mathcal{P}'_l$. Given a state digraph $\mathcal{D}(\bar{A}) = (\mathcal{X}, \mathcal{E}_{\bar{A}})$, we start with a minimal feasible solution $\mathcal{F}$. Notice that under such condition, a link failure may or may not compromise structural observability of $\mathcal{D}(\bar{A})$. Thus, we need to consider the properties of the links that jeopardizes structural observability when they fail. With these properties, we aim to identify new state variables that together with those from the starting minimal feasible solution form a $l$-robust feasible solution. Therefore, we introduce the notion of *sensitive link* with respect to $\mathcal{F}$.

**Definition 7.** *Consider a state digraph $\mathcal{D}(\bar{A}) = (\mathcal{X}, \mathcal{E}_{\bar{A}})$ and one of its minimal feasible solutions $\mathcal{F}$. We say $\mathcal{L}_{i,j}$ is a sensitive link with respect to $\mathcal{F}$ if and only if $\mathcal{F}$ is not a feasible solution of the corrupted state digraph $\mathcal{D}(\bar{A}_{-(i,j)})$.* $\qquad\square$

Hereafter, we can identify sensitive links by considering the minimum P&C decomposition, and the SCCs that compose the corrupted state digraph $\mathcal{D}(\bar{A}_{-(i,j)})$, as stated in the next result (that readily follows from Theorem 1).

**Lemma 2.** *Consider a state digraph $\mathcal{D}(\bar{A}) = (\mathcal{X}, \mathcal{E}_{\bar{A}})$ and one of its minimal feasible solutions $\mathcal{F}$. A link $\mathcal{L}_{i,j}$ is a sensitive link with respect to a minimal feasible solution $\mathcal{F}$ if and only if $\mathcal{F}$ does not contain the tips of the paths in a minimum P&C decomposition $\Xi$ of $\mathcal{D}(\bar{A}_{-(i,j)})$, or $\mathcal{F}$ does not contain at least one state variable in each sink-SCC of $\mathcal{D}(\bar{A}_{-(i,j)})$.* $\qquad\square$

To sharpen the intuition on sensitive links, consider the following example.

**Example 3.** *Consider the example in Fig. 2, which depicts a 5-state system. The state digraph consists of two SCCs, whose vertices are enclosed by the two dashed rectangles. The only sink-SCC contains $x_1$, and $x_5$. Fig. 2-(a) depicts one minimal feasible solution $\mathcal{F} = \mathcal{T} \cup \mathcal{S}^{\circ} = \{x_1, x_3\}$, where $\mathcal{T} = \{x_3\}$, $\mathcal{S}^{\circ} = \{x_1\}$, and after the failure of the link $\mathcal{L}_{5,1} = \{(x_5, x_1)\}$ or $\mathcal{L}_{4,2} = \{(x_4, x_2)\}$, $\mathcal{F} = \{x_1, x_3\}$ is not a feasible solution for the corrupted state digraph $\mathcal{D}(\bar{A}_{-(5,1)})$ or $\mathcal{D}(\bar{A}_{-(4,2)})$, since there does not exist a minimum P&C decomposition of $\mathcal{D}(\bar{A}_{-(4,2)})$ nor $\mathcal{D}(\bar{A}_{-(5,1)})$, with the tips of the paths belonging to the sink-SCC composed by $\{x_1, x_3\}$. Hence, the link $\mathcal{L}_{5,1}$ and $\mathcal{L}_{4,2}$ are two sensitive links. Similarly, Fig. 2-(b) depicts another minimal feasible solution $\tilde{\mathcal{F}} = \tilde{\mathcal{T}} \cup \tilde{\mathcal{S}}^{\circ} = \{x_4, x_5\}$, where $\tilde{\mathcal{T}} = \{x_4\}$, $\tilde{\mathcal{S}}^{\circ} = \{x_5\}$, and after the failure of the link $\mathcal{L}_{1,5} = \{(x_1, x_5)\}$ or $\mathcal{L}_{3,2} = \{(x_3, x_2)\}$, $\tilde{\mathcal{F}} = \{x_4, x_5\}$ is a non-feasible solution for the corrupted state digraph $\mathcal{D}(\bar{A}_{-(1,5)})$, or $\mathcal{D}(\bar{A}_{-(3,2)})$, since there does not exist a minimum P&C decomposition of $\mathcal{D}(\bar{A}_{-(1,5)})$ or $\mathcal{D}(\bar{A}_{-(3,2)})$, with the tips of the paths belonging to the SCC composed by $\{x_4, x_5\}$. Hence, the link $\mathcal{L}_{1,5}$ and $\mathcal{L}_{3,2}$ are two sensitive links.* $\qquad\square$

Therefore, as seen in Example 3, the set of sensitive links is based on a particular feasible solution, and there may exist several sensitive links. So, let the set of sensitive links associated with $\mathcal{F}$ be given by $\mathcal{L}(\mathcal{F}) = \cup_{j=1}^{\rho} \mathcal{L}_{l_j, r_j}$, where $\rho$ denotes the number of sensitive links.

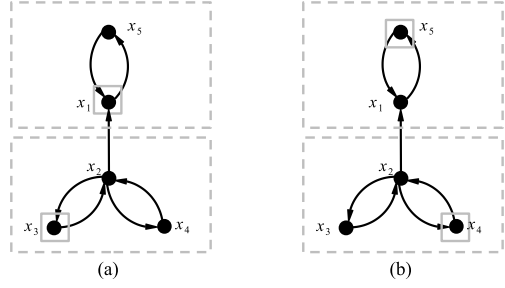

Fig. 2. A state digraph of a $5$-state system, where the grey squares depict the state vertices labeling state variables forming a solution, and the rectangles depict the two SCCs forming the state digraph. In (a), $\mathcal{F} = \{x_1, x_3\}$ is a minimal feasible solution, and $\mathcal{L}_{5,1} = \{(x_5, x_1)\}$ and $\mathcal{L}_{4,2} = \{(x_4, x_2)\}$ are two sensitive links. In (b), $\bar{\mathcal{F}} = \{x_4, x_5\}$ is another minimal feasible solution, and $\mathcal{L}_{1,5} = \{(x_1, x_5)\}$ and $\mathcal{L}_{3,2} = \{(x_3, x_2)\}$ are two sensitive links.

To solve $\mathcal{P}'_l$, we proceed similar to Section 4.1, given a feasible solution $\mathcal{F}$ we reduce the problem of constructing an $l$-robust feasible solution to a weighted set covering problem. We start by exploring the properties of sensitive links.

Recall that a minimal feasible solution $\mathcal{F}(\Xi)$ contains two subsets of state variables $\mathcal{T}(\Xi)$ and $\mathcal{S}^{\circ}(\Xi)$, associated with a P&C decomposion of $\mathcal{D}(\bar{A})$ as prescribed in Theorem 2. In fact, due to a sensitive link $\mathcal{L}_{l_j, r_j}$ failure, we have two mutually exclusive cases (as consequence of Lemma 2):

(i) there does not exist a minimum P&C decomposition $\Xi'$ of $\mathcal{D}(\bar{A}_{-(i,j)})$ such that $\mathcal{F}$ contains the tips of the paths of $\Xi'$. In particular, in this case, the sensitive link $\mathcal{L}_{l_j, r_j}$ necessarily belongs to a path in $\Xi$. Further, the number of paths associated with any minimum P&C decomposition $\Xi'$ increases by at most $1$ in comparison to $\Xi$.

(ii) there exists a minimum P&C decomposition $\Xi'$ of $\mathcal{D}(\bar{A}_{-(i,j)})$ such that $\mathcal{F}$ contains the tips of the paths of $\Xi'$, while the number of sink-SCCs without tips of the paths, associated with $\Xi$ increases by $1$.

We depict in Fig. 3-(b) and Fig. 3-(c) cases presented in (i) and (ii), respectively.

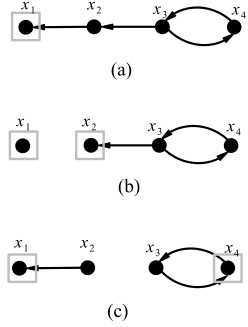

Fig. 3. An example of state digraphs, where the grey squares enclose the state vertices corresponding to a solution. In (a), $\mathcal{F} = \{x_1\}$ is a minimal feasible solution, where $\mathcal{T} = \{x_1\}$, $\mathcal{S}^\circ = \emptyset$. In (b), under the failure of $\mathcal{L}_{2,1}$, the number of tips of the paths for any minimum P&C decomposition is increased by $1$ with respect to $\mathcal{T}$, and $x_2$ is a additional tip of path which is required to belong to the $l$-robust feasible solution; finally, (c) depicts the scenario where a failure $\mathcal{L}_{3,2}$ occurs, where $\{x_1\}$ still constitutes the tips of the paths for some minimum P&C decomposition, but the SCC containing the vertices $x_3$ and $x_4$ requires one of its variables to be considered to obtain an $l$-robust feasible solution.

Based on the previous discussion, the set of sensitive links can be classified without loss of generality into two collections $\mathcal{L}_1$ and $\mathcal{L}_2$: $\mathcal{L}_1 = \cup_{j=1}^r \mathcal{L}_{l_j,r_j}$ represents the set of $r$ sensitive links corresponding to case (i), and $\mathcal{L}_2 = \cup_{j=r+1}^\rho \mathcal{L}_{l_j,r_j}$ represents the set of $\rho - r$ sensitive links corresponding to case (ii).

Subsequently, we introduce the notion of *completions* that correspond to a set of state variables $\mathcal{F}'$ such that after the failure of any sensitive link $\mathcal{L}_{l_j,r_j}$, $\mathcal{F} \cup \mathcal{F}'$ still forms a feasible solution for $\mathcal{D}(\bar{A}_{-(l_j,r_j)})$. Let us start by introducing the completions for elements in $\mathcal{L}_2$: after the failure of sensitive link $\mathcal{L}_{l_j,r_j} \subset \mathcal{L}_2$, $x_{l_j}$ and $x_{r_j}$ are separated into two different SCCs in $\mathcal{D}(\bar{A}_{-(l_j,r_j)})$. After such link failure, in $\mathcal{D}(\bar{A}_{-(l_j,r_j)})$, the SCC that $x_{l_j}$ belongs to becomes a sink-SCC without tip of paths. Hence, in order to obtain a feasible solution for $\mathcal{D}(\bar{A}_{-(l_j,r_j)})$, an additional state variable in the sink-SCC, that $x_{l_j}$ belongs to in $\mathcal{D}(\bar{A}_{-(l_j,r_j)})$, has to be considered. As a result, we can define the set of *sink-completions* $\Gamma_j^{\mathcal{S}}$ to a particular sensitive link $\mathcal{L}_{l_j,r_j} \subset \mathcal{L}_2$.

**Definition 8** (Sink-completions). *Consider a state digraph $\mathcal{D}(\bar{A}) = (\mathcal{X}, \mathcal{E}_{\bar{A}})$ and a minimal feasible solution $\mathcal{F} = \mathcal{T} \cup \mathcal{S}^{\circ}$. Further, consider the sensitive links $\mathcal{L}(\mathcal{F}) = \cup_{j=1}^{\rho} \mathcal{L}_{l_j, r_j} = \mathcal{L}_1 \cup \mathcal{L}_2$, where $\mathcal{L}_2 = \cup_{j=r+1}^{\rho} \mathcal{L}_{l_j, r_j}$ represents the set of sensitive links corresponding to case (ii). Then, for all $r < i \le \rho$,*

$$\Gamma_i^{\mathcal{S}} = \left\{ x' \in \mathcal{X} : x' \text{ belongs to the same SCC as } x_{l_i} \text{ in } \mathcal{D}(\bar{A}_{-(l_i, r_i)}) \right\}, \tag{13}$$

*is the set of sink-completions to sensitive link $\mathcal{L}_{l_i, r_i} \in \mathcal{L}_2$. Further, for notational convenience, we denote the state variables in $\Gamma_i^{\mathcal{S}}$ as $\Gamma_i^{\mathcal{S}} = \{\gamma_{\tau_i, 1}^{\mathcal{S}}, \cdots, \gamma_{\tau_i, n_i}^{\mathcal{S}}\}$, where $n_i$ denotes the total number of state variables in the sink-SCC of $\mathcal{D}(\bar{A}_{-(l_j, r_j)})$ that $x_{l_i}$ belongs to.* □

Similarly, we can introduce the notion of completions for sensitive links in $\mathcal{L}_1$, corresponding to case (i), which we refer to as *tip-completions*.

**Definition 9** (Tip-completions). *Consider a state digraph $\mathcal{D}(\bar{A}) = (\mathcal{X}, \mathcal{E}_{\bar{A}})$ and a minimal feasible solution $\mathcal{F} \equiv \mathcal{F}(\Xi) = \mathcal{T}(\Xi) \cup \mathcal{S}^{\circ}(\Xi)$. Further, consider the corresponding set of sensitive links $\mathcal{L}(\mathcal{F}) = \cup_{j=1}^{\rho} \mathcal{L}_{l_j, r_j} = \mathcal{L}_1 \cup \mathcal{L}_2$, where $\mathcal{L}_1 = \cup_{j=1}^{r} \mathcal{L}_{l_j, r_j}$ represents the set of sensitive links corresponding to case (i). Then, for all $i \le r$,*

$$\begin{aligned}\Gamma_i^{\mathcal{T}} = \{x \in \mathcal{X} : \; &\mathcal{T}(\Xi) \cup \{x\} \text{ contains a set of tips of the paths} \\ &\text{associated with a minimum P\&C decomposition of } \mathcal{D}(\bar{A}_{-(l_i, r_i)})\}\end{aligned} \tag{14}$$

*is the set of tip-completions to sensitive link $\mathcal{L}_{l_i, r_i} \in \mathcal{L}_1$. Further, for notational convenience, we denote the state variables in $\Gamma_i^{\mathcal{T}}$ as $\Gamma_i^{\mathcal{T}} = \{\gamma_{i,1}^{\mathcal{T}}, \cdots, \gamma_{i,\mu_i}^{\mathcal{T}}\}$, where $\mu_i$ is the number of tip-completions to $\mathcal{L}_{l_i, r_i}$.* □

Now, notice that given a minimal feasible solution $\mathcal{F}$, if we consider a sensitive link $\mathcal{L}_{l_i, r_i}$, by adding one corresponding completion $\gamma_i$, i.e., sink- or tip- completion, $\mathcal{F} \cup \{\gamma_i\}$ is a minimal feasible solution for the corrupted system $\bar{A}_{-(l_i, r_i)}$ after the failure of sensitive link $\mathcal{L}_{l_i, r_i}$.

**Remark 1.** *If the failure of sensitive link $\mathcal{L}_{l_i, r_i}$ of case (i) also increases the number of sink-SCCs that belong to the state digraph, then all elements of tip-completions $\Gamma_i^{\mathcal{T}}$ must belong to the new sink-SCC generated by the failure of $\mathcal{L}_{l_i, r_i}$.* □

 

In summary, we can generalize our previous discussion by considering collection of sets that account for the sink- and tip-completions for each sensitive link $\mathcal{L}_{l_i,r_i} \subset \mathcal{L}(\mathcal{F})$, which we refer to as *back-ups* and denoted by $\Theta^i_{\mathcal{F}}$. In particular, the sets in $\Theta^i_{\mathcal{F}}$ contain one state variable. Subsequently, we have the following definition.

**Definition 10** (Set of Back-ups). *Consider a state digraph $\mathcal{D}(\bar{A}) = (\mathcal{X}, \mathcal{E}_{\bar{A}})$ and a minimal feasible solution $\mathcal{F}$. Further, consider the set of $\rho$ sensitive links $\mathcal{L} = \cup^{\rho}_{j=1}\mathcal{L}_{l_j,r_j} = \mathcal{L}_1 \cup \mathcal{L}_2$, where $\mathcal{L}_1 = \cup^{r}_{j=1}\mathcal{L}_{l_j,r_j}$ represents the set of $r$ sensitive links corresponding to case (i), and $\mathcal{L}_2 = \cup^{\rho}_{j=r+1}\mathcal{L}_{l_j,r_j}$ represents the set of $\rho - r$ sensitive links corresponding to case (ii). Then, the set of back-ups to $\mathcal{L}_{l_i,r_i}$, referred to as $\Theta^i_{\mathcal{F}}$, is given by*

1) *for $i = 1, \ldots, r$,*

$$\Theta^i_{\mathcal{F}} = \left\{ \{\gamma^{\mathcal{T}}_{i,1}\}, \cdots, \{\gamma^{\mathcal{T}}_{i,\mu_i}\} \right\}, \tag{15}$$

   *where $\mu_i$ is the total number of tip-completions of $\mathcal{L}_{l_i,r_i}$ (see Definition 9); and*

2) *for $i = r+1, \ldots, \rho$,*

$$\Theta^i_{\mathcal{F}} = \left\{ \{\gamma^{\mathcal{S}}_{i,1}\}, \cdots, \{\gamma^{\mathcal{S}}_{i,n_i}\} \right\}, \tag{16}$$

   *where $n_i$ is the total number of sink-completions of $\mathcal{L}_{l_i,r_i}$ (see Definition 8).* □

In summary, $\Theta^i_{\mathcal{F}}$ is a collection of subsets of state variables each of which consists of one state variable, and (15)-(16) can be used to construct an $l$-robust feasible solution from a minimal feasible solution. More precisely, if we keep all state variables in minimal feasible solution $\mathcal{F}$ fixed, and for $\mathcal{D}(\bar{A}_{-(l_i,r_i)})$, we add any element in $\Theta^i_{\mathcal{F}}$ to the solution $\mathcal{F}$ (recall that it is no longer a feasible solution), then we achieve a feasible solution for $\mathcal{D}(\bar{A}_{-(l_i,r_i)})$. Consequently, if we have a set of state variables $\mathcal{F}'$, such that $\theta \subset \mathcal{F}'$ with $\theta \in \Theta^i_{\mathcal{F}}$, $\forall i \in \{1, \cdots, \rho\}$, then $\mathcal{F} \cup \mathcal{F}'$ is a feasible solution for any $\mathcal{D}(\bar{A}_{-(l_i,r_i)})$, i.e., $\mathcal{F}' \cup \mathcal{F}$ is an $l$-robust feasible solution.

Therefore, our goal is to find a minimal $\mathcal{F}'$ with the above properties, i.e., there is no other $\mathcal{F}''$ with $|\mathcal{F}''| < |\mathcal{F}'|$ such that $\mathcal{F} \cup \mathcal{F}''$ is an $l$-robust feasible solution. Similar to Section 4.1, a common back-up may exist for different elements in $\mathcal{L}$. As a result, minimizing the size of $\mathcal{F}'$ is the same as maximizing the number of shared state variables in back-up sets for different elements in the set of sensitive links $\mathcal{L}(\mathcal{F})$. Subsequently, we can find the minimal $\mathcal{F}'$ by considering an optimal set covering that covers the back-ups for any sensitive link in $\mathcal{L}(\mathcal{F})$. Note that in this section, the formulation of sets of back-ups is different from that of Section 4.1, where the elements of $\Theta^i_{\mathcal{F}}$ consist of a single state variable. Nevertheless, we can still use the similar notation as (8) to get a compressed expression of $\Theta^i_{\mathcal{F}}$. More specifically, for a system with $n$ state variables, we denote $\mathcal{Z}$ as $\mathcal{Z} = \{\mathcal{Z}_1, \cdots, \mathcal{Z}_n\}$, where

$$\mathcal{Z}_i = \{x_i\}, \text{ for } i = 1, \cdots, n.$$

Consequently, any element in $\Theta^i_{\mathcal{F}}$ consists of one element in $\mathcal{Z}$, which implies that $\Theta^i_{\mathcal{F}} = \{\mathcal{Z}_j\}_{j \in \mathcal{Q}}$, with $\mathcal{Q} \subset \{1, \ldots, n\}$. With this definition, the problem reduces to identifying for each element $\mathcal{Z}_i$, the subset of elements of $\mathcal{L}(\mathcal{F})$ for which $\mathcal{Z}_i$ is a shared back-up. This can be achieved by considering the definition of the sets $\mathcal{V}^l_j$ that contain the indices of elements in $\mathcal{L}(\mathcal{F})$, or equivalently $\Theta^i_{\mathcal{F}}$, that an element $\mathcal{Z}_j$ is a back-up to. Formally, this can be defined as follows: let $\mathcal{I} = \{1, \cdots, \rho\}$, then

$$\mathcal{V}^l_j = \{i \in \mathcal{I} : \mathcal{Z}_j \in \Theta^i_{\mathcal{F}}\}, \quad j = 1, \cdots, n, \tag{17}$$

i.e., the indices of the sets $\Theta^i_{\mathcal{F}}$'s to which the set $\mathcal{Z}_j$ belongs to.

As we can see in (15)-(16), the sets $\mathcal{Z}_i$'s are singletons, which is different from Section 4.1. As a result, in order to determine a minimal $l$-robust feasible solution $\mathcal{F}^l$ with $\mathcal{F}^l \supset S$, we can determine the minimum number of back-ups based on a particular feasible solution $\mathcal{F}$ that satisfies the robustness criterion. In other words, we cast $\mathcal{P}'_l$ as a set covering problem instead of a weighted set covering problem (nonetheless the problem is still NP-hard), as given in the following result.

**Theorem 7.** *Consider a state digraph $\mathcal{D}(\bar{A}) = (\mathcal{X}, \mathcal{E}_{\bar{A}})$ and one of its minimal feasible solutions $\mathcal{F} = \{x_{\tau_1}, x_{\tau_2}, \cdots, x_{\tau_p}\}$. Let the set of sensitive links associated with $\mathcal{F}$ be given by $\mathcal{L}(\mathcal{F}) = \cup^\rho_{j=1} \mathcal{L}_{l_j, r_j}$, where $\rho$ denotes the total number of sensitive links associated with $\mathcal{F}$. Consider the* *sets (17), and let $\mathcal{I} = \{1, \cdots, \rho\}$. If there exists $\mathcal{J}^*$ satisfying:*

$$\mathcal{J}^* = \underset{\mathcal{J} \subset \{1, \cdots, n\}}{\arg\min} \quad |\mathcal{J}|,$$

$$\text{subject to} \quad \mathcal{I} \subset \bigcup_{j \in \mathcal{J}^*} \mathcal{V}_j^l, \tag{18}$$

*i.e., the family $\{\mathcal{V}_j^l\}_{j \in \mathcal{J}^*}$ covers $\mathcal{I}$, then*

$$\mathcal{F}^* = \bigcup_{j \in \mathcal{J}^*} \mathcal{Z}_j \ \cup \ \mathcal{F},$$

*is an $l$-robust feasible solution. In addition, there is no other $l$-robust feasible solution $\mathcal{F}'$, with $\mathcal{F} \subset \mathcal{F}'$, such that $|\mathcal{F}'| < |\mathcal{F}^*|$.* □

The discussion in Section 4.1 about the computation efforts readily extends to Theorem 7, that requires to solve an NP-hard problem – the set covering problem. Consequently, we propose the solution to be approximated with polynomial complexity, which can be formally stated as follows.

**Theorem 8.** *Consider a state digraph $\mathcal{D}(\bar{A}) = (\mathcal{X}, \mathcal{E}_{\bar{A}})$, an arbitrary minimal feasible solution $\mathcal{F}$, the sets (17), and the set covering problem as described in Theorem 7. Then, the corresponding weighted set covering can be constructed with complexity $\mathcal{O}(|\mathcal{X}|^6)$. Further, a feasible but approximate solution $\mathcal{J}'$ of the corresponding minimum set covering may be constructed using an $\mathcal{O}(|\mathcal{X}|^3)$ complexity algorithm such that*

$$\mathcal{F}' = \bigcup_{j \in \mathcal{J}'} \mathcal{Z}_j \ \cup \ \mathcal{F},$$

*is an $l$-robust feasible solution, and the* performance ratio, *i.e., the ratio between the size of the approximated $l$-robust feasible solution and the size of minimal $l$-robust feasible solution containing $\mathcal{F}$ is bounded by the harmonic number $H(\rho) = \sum_{k=1}^{\rho} \frac{1}{k}$ of $\rho$, where $\rho$ is the total number of sensitive links.* □

Similar to the discussion in Section 4.1, we illustrate that the proposed two-step approach only ensures minimality of the $l$-robust feasible solution based on a particular minimal feasible solution. In other words, the minimality among all $l$-robust feasible solutions depends on the minimal feasible solution $\mathcal{F}$ that one starts with. Different minimal feasible solutions may lead to different sets of sensitive links, which results in different set covering problems, whose corresponding optimal solutions lead to $l$-robust feasible solutions with different sizes. Consider the example in Fig. 4, which depicts a 6-state system and two of its minimal feasible solutions with the dedicated outputs depicted by grey squares, respectively. Fig. 4-(a) depicts one minimal feasible solution $\mathcal{F} = \{x_1, x_6\}$, with no sensitive link. Hence, according to Theorem 7, $\mathcal{F} = \{x_1, x_6\}$ is also a minimal $l$-robust feasible solution. However, Fig. 4-(b) depicts another minimal feasible solution $\mathcal{F} = \{x_3, x_4\}$ where $\{(x_1, x_2)\}$ and $\{(x_6, x_5)\}$ are sensitive links. The sets of back-ups for the sensitive links are $\{\{x_1\}\}$ and $\{\{x_6\}\}$, respectively. Thereafter, $\mathcal{F} \cup \{x_1, x_6\}$ is not a minimal $l$-robust feasible solution.

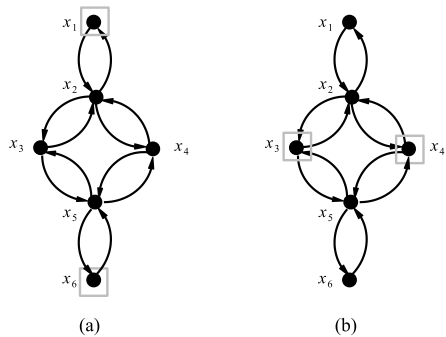

(a)  (b)

Fig. 4. This figure depicts a state digraph $\mathcal{D}(\bar{A})$ of a 6-state system, where the grey squares depict the state vertices labeling the state variables forming a solution. In (a) and (b) we show that $\mathcal{F} = \{x_1, x_6\}$ and $\mathcal{F}' = \{x_3, x_4\}$ are two minimal feasible solutions of $\mathcal{D}(\bar{A})$, which lead to $l$-robust feasible solutions with different sizes based on our approach.

Although our solution does not guarantee a minimal $l$-robust feasible solution, illustrative examples provided in Section 5 show that our set covering based design approach leads to $l$-robust feasible solutions with small number of states variables.

Next, we discuss how to extend the results on obtaining $l$-robust feasible solutions when directed links were considered to the case where we have undirected link failures.

### 4.2.2 Undirected Link Failure

Recall the two-step approach in Section 4.2.1 to compute an $l$-robust feasible solution when the links are directed. An undirected link consists of two directed links, so the procedure to find an $l$-robust feasible solution is similar to the directed one. More precisely, an *undirected link* between $x_j$ and $x_i$ is represented as $\mathcal{L}_{i,j} = \{(x_i, x_j), (x_j, x_i)\} \subset \mathcal{E}_{\bar{A}}$. The first step considers a minimal feasible solution used to determine the (undirected) sensitive links as described by Lemma 2. The set of undirected sensitive links still involves two cases: (i) the lack of the set of tips in a minimum P&C decomposition in $\mathcal{F}$ after an undirected link failure, and (ii) the generation of an additional sink-SCC without tips from the paths in the minimum P&C decomposition. Subsequently, we can enumerate completions for undirected sensitive links, and encode the enumeration of the completions as sets in a set covering problem as in Theorem 7. The sub-optimal solution with respect to undirected link failures can be obtained as part of the solution to the set covering problem, given by Theorem 8.

## 4.3 Extension to multiple sensor/link failures

As consequence of the results presented in Section 4.1 and Section 4.2, it follows that the problem of determining $s$-robust/$l$-robust feasible solutions with the minimum number of observed variables ensuring structural observability with respect to an arbitrary number of sensor/link failures is also NP- hard.

Further, we can extend the current results in Section 4.1 and Section 4.2 to address multiple sensor/link failures. Recall the two-step approach in Section 4.1 and Section 4.2, where the first step considers a minimal feasible solution, and in the second step, we find sets of back-ups for each state variable in the minimal feasible solution, or sets of completions for each sensitive link. Then, to obtain the minimum back-up and completion set we reduce the problem to a weighted set covering as we did in Section 4.1 and Section 4.2, respectively. Notwithstanding, when multiple failures occur, we need to find sets of back-ups (resp. completions) for each combination of $k$ state variables (resp. sensitive links). For instance, in the multiple sensor failure case, the universal set of the weighted set covering problem will be of size $\binom{p}{k}$, which corresponds to the total number of possible combinations of $k$ state variables in the minimal feasible solution, and the sets of back-ups contains subsets of state variables of size $k, k+1, \cdots, 2k$. More precisely, for each combination of $k$ state variables in the minimal feasible solution, at least we need $k$ sensors as a back-up, and at most $2k$ sensors is necessary as a back-up.

## 5 ILLUSTRATIVE EXAMPLE

In this section, we illustrate how to obtain an $s$-robust feasible solution and an $l$-robust feasible solution for a $10$-state system (randomly generated), whose state digraph is depicted in Fig. 5-(a). A minimal feasible solution is given by $\mathcal{F} = \{x_8, x_{10}\}$, whose corresponding state vertices are enclosed by the squares in Fig. 5-(b).

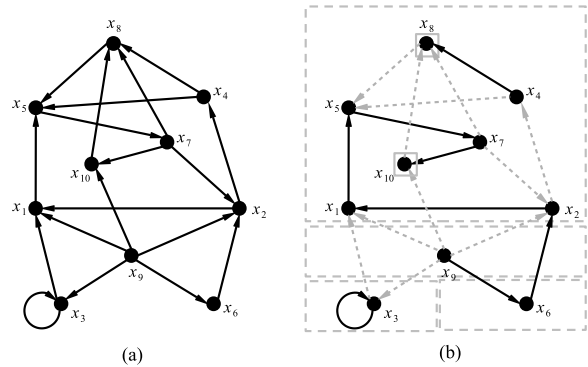

Fig. 5. In (a) we depict the state digraph with $10$ states. In (b) the four SCCs are depicted by the four dashed rectangles, and a minimal feasible solution is given by $\mathcal{F} = \{x_8, x_{10}\}$, whose corresponding state vertices are enclosed by the squares.

### 5.1 Robustness with respect to sensor failure

The back-ups sets (see (6)-(7)) are given by $\Omega_{\mathcal{F}} = \{\Omega_{\mathcal{F}}^1, \Omega_{\mathcal{F}}^2\}$, where $\Omega_{\mathcal{F}}^1 = \{\{x_1\}\}$; and $\Omega_{\mathcal{F}}^2 = \{\{x_1\}, \{x_4\}\}$. By invoking Theorem 4, we observe that there exist two tips of the paths, hence $\mathcal{I} = \{1, 2\}$, and the covering sets (see (9)-(10)) are given by $\mathcal{V}_1^s = \{1, 2\}$, $\mathcal{V}_4^s = \{1\}$, whereas the rest of the covering sets are empty, since all elements of $\Omega_{\mathcal{F}}^i$ consists of only one state variable. Subsequently, we can simply neglect those $\mathcal{V}_j^s$'s with $n < j \leq N$, where $n = 10$ is the total number of state variables, and $N = n + \binom{n}{2} = \frac{n(n+1)}{2}$ is the total number of sets that consist of one or two state variables.

As a result, we constructed a set covering problem as proposed in Theorem 4, with the universal set $\mathcal{I}$, and the covering sets $\mathcal{V}_i^s$, $i \in \{1, \cdots, 10\}$. In this example, we can readily identify that $\mathcal{V}_1^s$ is an optimal solution to the set covering problem in Theorem 4, which implies that $\mathcal{F}^* = \mathcal{F} \cup \{x_1\} = \{x_1, x_8, x_{10}\}$ is an $s$-robust feasible solution. However, for larger systems, we may need to consider sub-optimal solutions as discussed in Theorem 5.

## 5.2   Robustness with respect to link failure

First, we identify the set of directed sensitive links associated with the choice of $\mathcal{F}$. This consists of $5$ elements, given by $\mathcal{L}_{1,5}^{(d)}, \mathcal{L}_{5,7}^{(d)}, \mathcal{L}_{6,2}^{(d)}, \mathcal{L}_{7,10}^{(d)}$, and $\mathcal{L}_{3,1}^{(d)}$.

For those sensitive links, the collection of back-ups sets (see (15)-(16)) is given by $\Theta_{\mathcal{F}} = \{\Theta_{\mathcal{F}}^1, \Theta_{\mathcal{F}}^2, \cdots, \Theta_{\mathcal{F}}^5\}$, where $\Theta_{\mathcal{F}}^1 = \{\{x_1\}\}$, $\Theta_{\mathcal{F}}^2 = \{\{x_5\}\}$, $\Theta_{\mathcal{F}}^3 = \{\{x_6\}\}$, $\Theta_{\mathcal{F}}^4 = \{\{x_6, x_7\}\}$, and $\Theta_{\mathcal{F}}^5 = \{\{x_3\}\}$.

Because there exist five sensitive links, we have $\mathcal{I} = \{1, \cdots, 5\}$, and the covering sets (see (17)) are given by $\mathcal{V}_1^l = \{1\}$, $\mathcal{V}_3^l = \{6\}$, $\mathcal{V}_5^l = \{2\}$, $\mathcal{V}_6^l = \{3, 4\}$, $\mathcal{V}_7^l = \{4\}$, and the other sets are empty. As a result, we constructed a set covering problem as proposed in Theorem 7, with the universal set $\mathcal{I}$, and the covering sets $\mathcal{V}_i^l$, for $i \in \{1, \cdots, 10\}$. Further, we can readily identify that either $\{\mathcal{V}_i^l\}_{i \in \{1,3,5,6\}}$ is an optimal solution to the set covering problem in Theorem 7. Thus, $\mathcal{F}^* = \mathcal{F} \cup \{x_1\} = \{x_1, x_3, x_5, x_6, x_8, x_{10}\}$ is an $l$-robust feasible solution. However, for larger systems, we may need to consider sub-optimal solutions as discussed in Theorem 8.

## 5.3   Discussion of Results

To evaluate the trade-off between the running time of the proposed algorithm and the network topology, we consider state digraphs that were randomly generated using three models commonly found in the literature [32]: Erdős-Rényi [33], small world and scale free models.

In our simulation, we consider two properties of random networks: the scale of the network and the connectivity of the network. Simulation results show that in Erdős-Rényi model, usually we do not have any $s$-robust feasible solution when the connectivity of the network is small, in which case the graph contains several source SCCs with only one state vertex, which implies that the corresponding state variable has to be measured directly. Alternatively, when state digraphs are modeled as small world networks, due to the likelihood of existing a cycle which connects all state vertices, all minimal feasible solution consist of one sensor, and we need at most two sensors to ensure $s$-robustness or $l$-robustness. As a result, we now devil in details when the state digraphs have a scale-free network topology.

In the scale-free network implemented as in [34], there are two parameters: $n$ denotes the number of nodes in graph, and $d$ denotes the minimum node degree, which determines the connectivity of the network. In addition, we randomly switched $10\%$ of undirected links to directed links. For the simulation experiments, we used a Macbook Pro running Ubuntu Linux with a 2.7 GHz Intel Core i5 processor. In order to compute a minimum P&C decomposition, toolbox TOMLAB/CPLEX [35] was considered. After we construct the set covering problems as in Theorem 4 and Theorem 7, we use greedy algorithm [29] to find a sub-optimal solution to the set covering problems. The Matlab implementation of the proposed algorithm can be found in [36].

In Figure 6, we show the simulation results when we consider different values of $d$ with $n = 300$, and record the runtime and the cardinalities of solutions based on proposed algorithms. In Figure 7, we show the simulation results when we consider different values of $n$ with $d = 1$, and record the runtime and the cardinalities of solutions based on proposed algorithms. Furthermore, we fitted the curves in Figure 7-(a) with functions of the form $f(n) = a * n^b$. For $s$-robust, the values of $a, b$ with $b \in \mathbb{N}$ that lead to the minimal square of the correlation are $a = 2.235 \times 10^{-7}$ and $b = 4$. For $l$-robust, the values of $a, b$ with $b \in \mathbb{N}$ that lead to the minimal square of the correlation are $a = 2.673 \times 10^{-7}$ and $b = 4$.

In addition, Table I lists the simulation results, where we consider different values of $n$ with $d = 1$, and record the number of minimum P&C decomposition required for

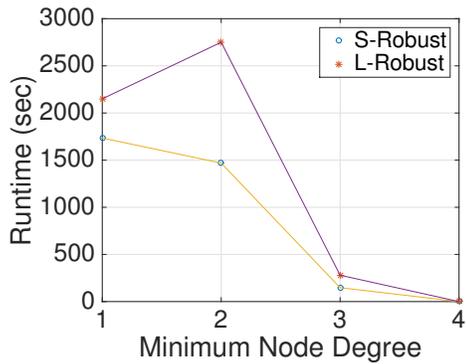
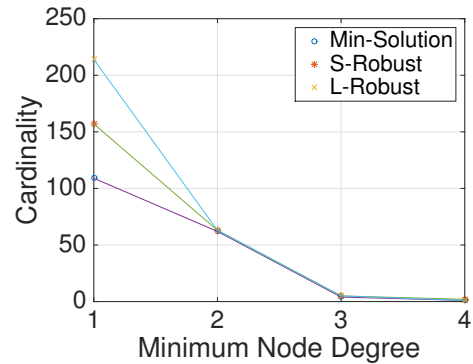

(a) Runtime to compute $s/l$-robust feasible solutions

(b) Cardinality of minimal and $s/l$-robust feasible solutions

Fig. 6. This figure depicts the relationship between the minimum node degree and the runtime/cardinalities of solutions for scale-free networks with different values of $d$ and $n = 300$. In (a), the relationship between the minimum node degree and the runtime is considered. In (b), the relationship between the minimum node degree and the cardinalities of solutions is depicted.

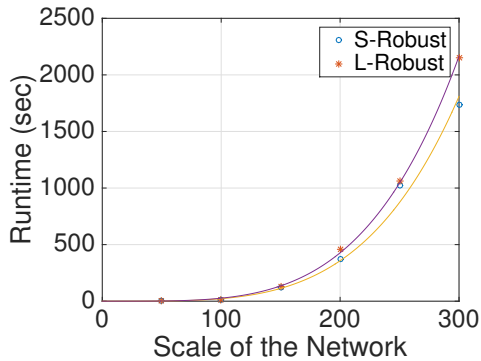
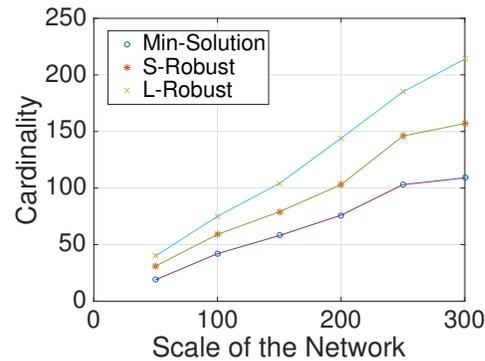

(a) Runtime to compute $s/l$-robust feasible solutions

(b) Cardinality of minimal and $s/l$-robust feasible solutions

Fig. 7. This figure depicts the relationship between the minimum node degree and the runtime/cardinalities of solutions for scale-free networks with different values of $n$ and $d = 1$. In (a), the relationship between the scale of the network and the runtime is considered. In (b), the relationship between the scale of the network and the cardinalities of solutions is depicted.

each experiment. In Table I, $D_s$ represents the number of P&C decomposition required to construct the set covering problem in Theorem 4, $D'_s$ represents the number of P&C decomposition required to construct the set covering problem based on pure enumeration for $s$-robustness, $D_l$ represents the number of P&C decomposition required to construct the set covering problem in Theorem 7, $D'_l$ represents the number of P&C decomposition required to construct the set covering problem based on pure enumeration for $l$-robustness. Thus, from Table I it readily follows that the proposed approach enables a systematic and efficient method to determine $s/l$-robust feasible solutions. In addition, from Table I and Figure 7 we can see that the network with $n = 50$ already lead to a set covering problem with the universal set of size 19, and 558 covering sets. To find the optimal solution, in the worst case, combination of covering sets need to be computed is $\binom{558}{1} + \binom{558}{2} + \ldots + \binom{558}{19}$, which is cumbersome to compute. Hence, approximation algorithms like greedy algorithm are required to find an $s/l$ robust solution. Finally, although we present simulations only for models randomly generated, the proposed approach is general and applicable to those systems not modeled by the mentioned models.

| $n$ | 50 | 100 | 200 | 300 |
|---|---|---|---|---|
| $D_s$ | 558 | 2378 | 9300 | 20628 |
| $D'_s$ | $2.15 \times 10^9$ | $2.88 \times 10^{17}$ | $2.13 \times 10^{37}$ | $3.14 \times 10^{57}$ |
| $D_l$ | 651 | 2320 | 10416 | 26740 |
| $D'_l$ | $2.15 \times 10^9$ | $2.88 \times 10^{17}$ | $2.13 \times 10^{37}$ | $3.14 \times 10^{57}$ |

TABLE 1

# 6 CONCLUSIONS AND FURTHER RESEARCH

In this paper, we addressed the problem of deploying the minimum number of dedicated sensors that ensure a network of interconnected linear dynamical systems to be structurally observable under disruptive scenarios. The disruptive scenarios considered are: (i) the malfunction/loss of one arbitrary sensor, and (ii) the failure of connection (either directed or undirected) between a pair of dynamical devices in an interconnected dynamical system. Due to the combinatorial nature of the problem, we proposed a two step approach. First, we achieve an arbitrary minimum sensor placement configuration ensuring structural observability. Secondly, we develop a sequential process to find the minimum number of additional sensors required for robust observability with respect to arbitrary sensor/link failure. This step can be solved by recasting it as a (weighted) set covering problem, for which sub-optimal solutions can be determined in polynomial-time and with optimality guarantees.

Future research consists of finding efficient approaches to determine sensor deployment to cope with disruptive scenarios where multiple failures can occur at the same time.

# 7 APPENDIX

*Proof of Theorem 3:*

To prove that $\mathcal{P}'_s$ is NP-hard, we need to provide a reduction from an NP-hard problem to $\mathcal{P}'_s$. Intuitively, it would mean that if a solution to general instances of $\mathcal{P}'_s$ can be efficiently obtained, we would also obtain an efficient (polynomial) solution to any instance of the NP-hard problem. Therefore, consider the set covering problem with universal set $\mathcal{U} = \{1, \cdots, p\}$ and a finite collection of $k$ sets $\mathcal{C}_i \subset \mathcal{U}$ with $i \in \mathcal{J} = \{1, \cdots, k\}$. In addition, consider the reduction proposed in Fig. 8, where an instance of the set covering problem with input $(\mathcal{C}_1, \ldots, \mathcal{C}_k; \mathcal{U})$ is written as an instance to $\mathcal{P}'_s$, with input of the $(p+k+4) \times (p+k+4)$ matrix $\bar{A}(\mathcal{C}_1, \ldots, \mathcal{C}_k; \mathcal{U})$, whose state digraph $\mathcal{D}(\bar{A})$ is depicted in Figure 8, given as follows: (1) $\bar{A}_{p+j,i} = 1$ if $i \in \mathcal{C}_j$; (2) $\bar{A}_{p+j,p+j} = 1$ for $j = 1, \ldots, k$; (3) $\bar{A}_{p+k+2,p+j} = 1$ for $j = 1, \ldots, k$; (4) $\bar{A}_{p+k+2,p+k+1} = \bar{A}_{p+k+1,p+k+2} = \bar{A}_{p+k+3,p+k+2} = \bar{A}_{p+k+2,p+k+3} = \bar{A}_{p+k+2,p+k+4} = 1$; and zero otherwise.

Notice that $\mathcal{F}^s$ is an $s$-robust feasible solution only if it contains two state variables from the sink-SCC, i.e., two state variables out of $\{x_{p+k+1}, x_{p+k+2}, x_{p+k+3}\}$, recall Theorem 1. In addition, $\mathcal{F}^s \setminus \{x\}$, with $x \in \mathcal{F}^s$, has to contain the tips of the paths $\mathcal{T}(\Xi)$ associated with some minimum P&C decomposition $\Xi$. Due to the sparse nature of $\bar{A}(\mathcal{C}_1, \ldots, \mathcal{C}_k; \mathcal{U})$, it follows that the tips of the paths associated with a mini-mum P&C decomposition can be characterized in closed-form. More precisely, let $\mathcal{T}^1 = \{x_1, \ldots, x_p, x_{p+k+1}, x_{p+k+4}\}$ and $\mathcal{T}^2 = \{x_1, \ldots, x_p, x_{p+k+3}, x_{p+k+4}\}$ be two sets of tips of the paths and $\mathcal{R}_i = \{x_j : x_i \text{ reaches } x_j\}$ be the set of state variables reachable from $x_i$ (including itself). Then any set of tips of the paths is given by $\mathcal{T}^1 \setminus \{x_l\} \cup \{r\}$, where $x_l \in \mathcal{T}^1$ and $r \in \mathcal{R}_l \setminus (\mathcal{T}^1 \cup \{x_{p+k+2}\})$ or $\mathcal{T}^2 \setminus \{x_l\} \cup \{r\}$, where $x_l \in \mathcal{T}^2$ and $r \in \mathcal{R}_l \setminus (\mathcal{T}^2 \cup \{x_{p+k+2}\})$. Consequently, a minimum $s$-robust feasible solution has to contain $x_{p+k+1}$ and $x_{p+k+3}$ to ensure that two state variables from the sink-SCC are chosen, and such that they can account for one state variable in any set of tips of the paths; in other words, if the choice included $x_{p+k+2}$, then an additional state variable would be required to be included in an $s$-robust feasible solution. Further, from the description of the sets of tips of the paths, it follows that $\mathcal{F}^s$ has to contain the minimum collection of state variables in $\mathcal{R}_l \setminus (\mathcal{T} \cup \{x_{p+k+2}\})$ for all $l \in \{1, \ldots, p\}$. In addition, from the interpretation in terms of the set covering problem presented in Fig. 8, it follows that a solution to our problem would allow us to retrieve a solution to the set covering problem $(\mathcal{C}_1, \ldots, \mathcal{C}_k; \mathcal{U})$, and the result follows. $\diamond$

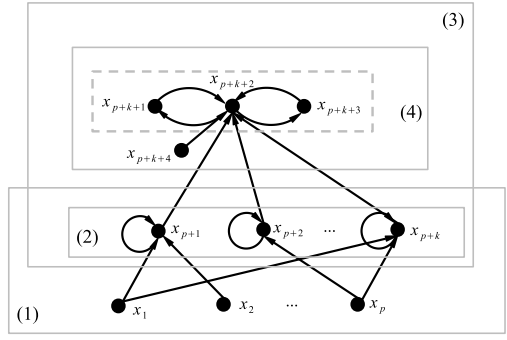

Fig. 8. This figure depicts the state digraph $\mathcal{D}(\bar{A}(\mathcal{C}_1, \ldots, \mathcal{C}_k; \mathcal{U}))$ with one sink-SCC with its vertices enclosed in the dashed box. Briefly, $x_i$ ($i = 1, \ldots, p$) are associated with the elements in the universal set, and $x_{p+j}$ ($j = 1, \ldots, k$) with the sets in the set covering problem. In other words, there is an edges from $x_i$ ($i = 1, \ldots, p$) to $x_{p+j}$ ($j = 1, \ldots, k$) if an element $i$ belongs to set $\mathcal{C}_j$.

*Proof of Lemma 1:* It suffices to verify that $\mathcal{F}^* \setminus \{x\}$ is a feasible solution for all $x \in \mathcal{F}^*$. To this end, we consider the following two cases: (i) if $x \in \mathcal{F}^* \setminus \mathcal{F}$, then $\mathcal{F} \subset (\mathcal{F}^* \setminus \{x\})$, and because $\mathcal{F}$ is a feasible solution, $\mathcal{F}^* \setminus \{x\}$ is also a feasible solution; (ii) if $x \in \mathcal{F}$, say $x = x_{\tau_i}$, it follows that $\bigcup_{j \in \mathcal{J}} \mathcal{Z}_j \subset (\mathcal{F}^* \setminus \{x\})$. Further, by definition of set covering and $\mathcal{V}_j^{s'}$'s in (9)-(10), there exists $j^* \in \mathcal{J}$, such that $i \in \mathcal{V}_{j^*}^s$, and $\mathcal{Z}_{j^*} \in \Omega_{\mathcal{F}}^i$. In other words, there exists a set of state variable(s) $\mathcal{Z}_{j^*} \subset (\mathcal{F}^* \setminus \{x\})$, such that $\{\mathcal{F} \setminus \{x_{\tau_i}\}\} \bigcup \mathcal{Z}_{j^*}$ is a feasible solution (by definition of $\Omega_{\mathcal{F}}^i$ in (6)-(7)). Hence, $\mathcal{F}^* \setminus \{x\}$ is a feasible solution. ■

*Proof of Theorem 4:* To prove optimality, we proceed as follows: (i) we show that $\mathcal{F}^*$ is an $s$-robust feasible solution; and (ii) we show that there is no other $s$-robust feasible solution $\mathcal{F}'$, with $\mathcal{F} \subset \mathcal{F}'$, such that $|\mathcal{F}'| < |\mathcal{F}^*|$. Condition (i) holds by Lemma 1. To prove (ii), we assume by contradiction that there exists an $s$-robust feasible solution $\mathcal{F}'$, with $\mathcal{F} \subset \mathcal{F}'$, such that $|\mathcal{F}'| < |\mathcal{F}^*|$. By definition of the sets $\mathcal{V}_j^s$ and the associated cost $c_j$ in (9)-(11), there exists a set $\mathcal{J}' \subset \{1, \cdots, N\}$ associated with $\mathcal{F}' \setminus \mathcal{F}$ such that $\mathcal{I} \subset \bigcup_{j \in \mathcal{J}'} \mathcal{V}_j^s$, and $\sum_{k \in \mathcal{J}'} c_k < \sum_{k \in \mathcal{J}^*} c_k$. This contradicts the fact that $\mathcal{J}^* = \arg\min_{\mathcal{J} \subset \{1, \cdots, N\}} \sum_{k \in \mathcal{J}} c_k$. As a result, $\bigcup_{j \in \mathcal{J}^*} \mathcal{Z}_j$ is the minimal set of back-ups for $\mathcal{F}$, and the result follows. ■

*Proof of Theorem 5:* A (minimal) feasible solution $\mathcal{F}$ can be efficiently determined using a polynomial complexity algorithm [6], with a time complexity of $\mathcal{O}(|\mathcal{X}|^3)$. Similarly, computing each set of $\Omega_{\mathcal{F}}^i$ for $i = 1, \cdots, m$ has a time complexity of $\mathcal{O}(|\mathcal{X}|^4)$. As a result, the complexity of computing $\Omega_{\mathcal{F}}^i$ for all $i = 1, \cdots, m$ is $\mathcal{O}(|\mathcal{X}|^5)$. Also, each set of $\Omega_{\mathcal{F}}^i$ for $i = m + 1, \cdots, p$ can be determined with a time complexity of $O(|\mathcal{X}|)$. Hence, the complexity of determining $\Omega_{\mathcal{F}}^i$ for all $i = m + 1, \cdots, p$ is $\mathcal{O}(|\mathcal{X}|^2)$. Remark that $\mathcal{V}_i^s$ can be efficiently implemented in at most $\mathcal{O}(|\mathcal{X}|^3)$ since it consists of verifying if each of the $|\mathcal{Z}|$ sets $\left(\frac{|\mathcal{X}|(|\mathcal{X}|+1)}{2}\right.$ in total$\left.\right)$ belongs to at most $|\mathcal{X}|$ sets of $\Omega_{\mathcal{F}}^i$'s. Thus, we can construct the corresponding weighted set covering with complexity $\mathcal{O}(|\mathcal{X}|^5)$. Further, for the weighted set covering problem we constructed, the size of the universal set is bounded by $|\mathcal{X}|$, and the number of cover sets is bounded by $|\mathcal{X}|^2$, which leads to our conclusion about the overall computational complexity.

By considering the greedy heuristic algorithm for the weighted set covering problem [37], we obtain a performance ratio of $H(d)$, where $d$ is the size of the largest cover set. Thus, the solution obtained with the proposed method achieves a performance ratio of $H(p)$ where $p$ is the size of the universal set. In addition, we use the fact that the algorithm in [37] has computation complexity bounded by $\mathcal{O}(p \times N)$, where $N$ represents the number of cover sets. Finally, we notice that feasibility holds due to Lemma 1. ■

*Proof of Theorem 6:*

We provide a similar reduction to that in Theorem 3, where an instance of the set covering problem with input $(\mathcal{C}_1, \ldots, \mathcal{C}_k; \mathcal{U})$ is written as an instance to $\mathcal{P}'_l$, with input of the $(p+2k+4) \times (p+2k+4)$ matrix $\bar{A}(\mathcal{C}_1, \ldots, \mathcal{C}_k; \mathcal{U})$, whose state digraph $\mathcal{D}(\bar{A})$ is depicted in Figure 9, given as follows: (1) $\bar{A}_{i,i} = 1$ for $i = 1, \ldots, p$ (2) $\bar{A}_{p+j,i} = 1$, and $\bar{A}_{p+k+j,i} = 1$ if $i \in \mathcal{C}_j$; (3) $\bar{A}_{p+j,p+j} = 1$ and $\bar{A}_{p+k+j,p+k+j} = 1$ for $j = 1, \ldots, k$; (4) $\bar{A}_{p+2k+2,p+j} = 1$ and $\bar{A}_{p+2k+2,p+k+j} = 1$ for $j = 1, \ldots, k$; (5) $\bar{A}_{p+2k+1,p+2k+3} = \bar{A}_{p+2k+1,p+2k+4} = \bar{A}_{p+2k+2,p+2k+3} = \bar{A}_{p+2k+2,p+2k+4} = 1$; and zero otherwise.

Now, notice that $\mathcal{F}^l$ is an $l$-robust feasible solution only if it contains state variables $\mathcal{T} \equiv \{x_{p+2k+1}, x_{p+2k+2}\}$ that is the tips of the paths associated with any minimum P&C decomposition of the state digraph. Further, $\mathcal{T}$ consists of a state variable from each sink-SCC. As a result, $\mathcal{T}$ is the only minimal feasible solution of $\bar{A}(\mathcal{C}_1, \ldots, \mathcal{C}_k; \mathcal{U})$. In addition, only links $\{(x_i, x_i) : i \in \mathcal{U}\}$ are sensitive links with respect to $\mathcal{T}$; in particular, $\mathcal{T}' = \mathcal{T} \cup \{x_i\}$ is the tips of the paths associated with a minimum P&C decomposition associated with the *corrupted* state digraph $\mathcal{D}(\bar{A}_{-(i,j)})$. Due to the sparse nature of $\bar{A}(\mathcal{C}_1, \ldots, \mathcal{C}_k; \mathcal{U})$, we can characterize all possible sets of tips of the paths associated with a minimum P&C decomposition of the corrupted state digraph. More precisely, let $\mathcal{R}_i = \{x_j : \text{ exists a direct path from } x_i \text{ to } x_j\}$ be the set of state variables whose corresponding vertices are reachable from $x_i$ (including itself) in $\mathcal{D}(\bar{A}(\mathcal{C}_1, \ldots, \mathcal{C}_k; \mathcal{U}))$, then any set of tips of the paths associated with the corrupted state digraph after the failure of sensitive link $\{(x_i, x_i)\}$ is given by $\{x_{p+2k+1}, x_{p+2k+2}\} \cup \{r\}$, where $r \in \mathcal{R}_i \setminus \{x_{p+2k+2}\}$. Consequently, from the description of the sets of tips of the paths, it follows that $\mathcal{F}^l$ has to contain the minimum collection of state variables in $\mathcal{R}_i$ for all $i \in \{1, \ldots, p\}$, where due to symmetry choosing $x_{p+i}$ or $x_{p+k+i}$ leads to the same results. In addition, from the interpretation in terms of the set covering problem presented in Fig. 9, it follows that a solution to our problem would allow us to retrieve a solution to the set covering problem $(\mathcal{C}_1, \ldots, \mathcal{C}_k; \mathcal{U})$, and the result follows. ∎

*Proof of Theorem 7:* Similar to Theorem 4, by considering the sensitive links and the sets (15)-(16) for the minimum set covering problem. ∎

*Proof of Theorem 8:* Similar to Theorem 5, by considering the sensitive links and the

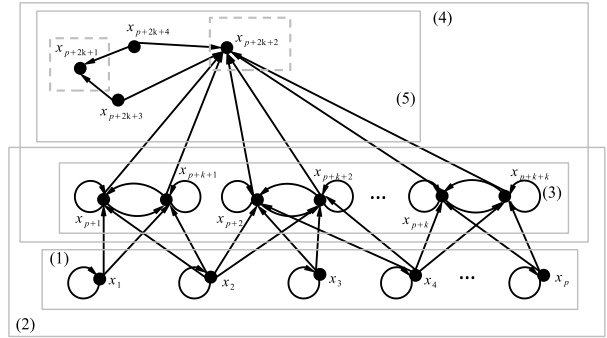

Fig. 9. This figure depicts the state digraph $\mathcal{D}(\bar{A}(\mathcal{C}_1, \ldots, \mathcal{C}_k; \mathcal{U}))$ with two sink-SCCs with its vertices enclosed in the two dashed boxes. Briefly, $x_i$ ($i = 1, \ldots, p$) are associated with the elements in the universal set, and $x_{p+j}, x_{p+k+j}$ ($j = 1, \ldots, k$) with the sets in the set covering problem. In other words, there are edges from $x_i$ ($i = 1, \ldots, p$) to $x_{p+j}$ and $x_{p+k+j}$ ($j = 1, \ldots, k$) if an element $i$ belongs to set $\mathcal{C}_j$.

sets (15)-(16) for the minimum set covering problem. ∎

## REFERENCES

[1] Y. Mo, T.-H. Kim, K. Brancik, D. Dickinson, L. Heejo, A. Perrig, and B. Sinopoli, "Cyber-physical security of a smart grid infrastructure," *Proceedings of the IEEE*, vol. 100, no. 1, pp. 195–209, Jan 2012.

[2] A. A. Cárdenas, S. Amin, and S. Sastry, "Research challenges for the security of control systems," in *Proceedings of the 3rd Conference on Hot Topics in Security*. Berkeley, CA, USA: USENIX Association, 2008.

[3] F. Pasqualetti, F. Dorfler, and F. Bullo, "Attack detection and identification in cyber-physical systems," *IEEE Transactions on Automatic Control*, vol. 58, no. 11, pp. 2715–2729, Nov 2013.

[4] A. Atputharajah and T. Saha, "Power system blackouts - literature review," in *2009 International Conference on Industrial and Information Systems (ICIIS)*, Dec 2009, pp. 460–465.

[5] A. Olshevsky, "Minimal controllability problems," *IEEE Transactions on Control of Network Systems*, vol. 1, no. 3, pp. 249–258, Sept 2014.

[6] S. Pequito, S. Kar, and A. Aguiar, "A framework for structural input/output and control configuration selection in large-scale systems," *IEEE Transactions on Automatic Control*, vol. PP, no. 99, pp. 1–1, 2015.

[7] A. Olshevsky, "Minimum input selection for structural controllability," in *In the Proceedings of the American Control Conference (ACC)*, July 2015, pp. 2218–2223.

[8] Y. Liu, J. Slotine, and A. Barabasi, "Observability of complex systems," *PNAS*, vol. 110(7), pp. 2460–2465, 2013.

[9] S. Pequito, G. Ramos, S. Kar, A. P. Aguiar, and J. Ramos, "On the Exact Solution of the Minimal Controllability Problem," *ArXiv e-prints*, 2014.

[10] T. Summers, F. Cortesi, and J. Lygeros, "On submodularity and controllability in complex dynamical networks," *IEEE Transactions on Control of Network Systems*, vol. PP, no. 99, pp. 1–1, 2015.

[11] V. Tzoumas, M. Rahimian, G. Pappas, and A. Jadbabaie, "Minimal actuator placement with bounds on control effort," *IEEE Transactions on Control of Network Systems*, vol. PP, no. 99, pp. 1–1, 2015.

[12] J.-M. Dion, C. Commault, and J. V. der Woude, "Generic properties and control of linear structured systems: a survey," *Automatica*, vol. 39, no. 7, pp. 1125 – 1144, 2003.

[13] S. Pequito, S. Kar, and A. P. Aguiar, "Minimum cost input/output design for large-scale linear structural systems," *Automatica*, vol. 68, pp. 384 – 391, 2016.

[14] C. Commault and J.-M. Dion, "Input addition and leader selection for the controllability of graph-based systems," *Automatica*, vol. 49, no. 11, pp. 3322 – 3328, 2013.

[15] S. Pequito, S. Kar, and A. P. Aguiar, "On the complexity of the constrained input selection problem for structural linear systems," *Automatica*, vol. 62, pp. 193 – 199, 2015.

[16] C. Rech, "Robustness of interconnected systems to structural disturbances in structural controllability and observability," *International Journal of Control*, vol. 51, no. 1, pp. 205–217, 1990.

[17] T. H. Do, C. Commault, and J. M. Dion, "Sensor classification for observability of structured systems," in *Proceedings of 7th Workshop on Advanced Control and Diagnosis*, 2009.

[18] C. Commault and J. M. Dion, "Sensor location for diagnosis in linear systems: A structural analysis," *IEEE Transactions on Automatic Control*, vol. 52, no. 2, pp. 155–169, 2007.

[19] V. Pichai, M. E. Sezer, and D. D. Šiljak, "Vulnerability of dynamic systems," *International Journal of Control*, vol. 34, no. 6, pp. 1049–1060, 1981.

[20] B. Wang, L. Gao, Y. Gao, and Y. Deng, "Maintain the structural controllability under malicious attacks on directed networks," *EPL (Europhysics Letters)*, vol. 101, no. 5, p. 58003, 2013.

[21] J. Ruths and D. Ruths, "Robustness of network controllability under edge removal," in *Complex Networks IV*, ser. Studies in Computational Intelligence, G. Ghoshal, J. Poncela-Casasnovas, and R. Tolksdorf, Eds. Springer Berlin Heidelberg, 2013, vol. 476, pp. 185–193.

[22] M. A. Rahimian and A. G. Aghdam, "Structural controllability of multi-agent networks: Robustness against simultaneous failures," *Automatica*, vol. 49, no. 11, pp. 3149 – 3157, 2013.

[23] X. Liu, S. Pequito, S. Kar, Y. Mo, B. Sinopoli, and A. Aguiar, "Minimum robust sensor placement for large scale linear time-invariant systems: a structured systems approach," *4th IFAC Workshop on Distributed Estimation and Control in Networked Systems (NecSys)*, pp. 417–424, 2013.

[24] F. Dorfler and M. R. Jovanovic, "Wide area control of ieee 39 new england power grid model." [Online]. Available: www.umn.edu/ mihailo/software/lqrsp/

[25] M. Mesbahi and M. Egerstedt, *Graph Theoretic Methods in Multiagent Networks*. Princeton University Press, 2010.

[26] D. D. Siljak, *Large-Scale Dynamic Systems: Stability and Structure*. Dover Publications, 2007.

[27] J. M. Dion, C. Commault, and J. V. der Woude, "Characterization of generic properties of linear structured systems for efficient computations," *Kybernetika*, vol. 38, pp. 503–520, 2002.

[28] T. H. Cormen, C. Stein, R. L. Rivest, and C. E. Leiserson, *Introduction to Algorithms*, 2nd ed. McGraw-Hill Higher Education, 2001.

[29] V. Vazirani, *Approximation Algorithms*. Springer-Verlag, 2001.

[30] D. S. Hochbaum, "Approximation algorithms for the set covering and vertex cover problems," *SIAM Journal on Computing*, vol. 11(3), pp. 555–556, 1982.

[31] M. A. Astrand and J. Suomela, "Fast distributed approximation algorithms for vertex cover and set cover in

anonymous networks," *Proceedings of the 22nd ACM symposium on Parallelism in algorithms and architectures*, pp. 294–302, 2010.

[32] R. van der Hofstad, "Random graphs and complex networks," 2011. [Online]. Available: www.win.tue.nl/rhofstad/NotesRGCN2011.pdf

[33] P. Erdős and A. Rényi, "On random graphs. I," *Publicationes Mathematicae*, vol. 6, pp. 290–297, 1959.

[34] V. Batagelj and U. Brandes, "Efficient generation of large random networks," *Phys. Rev. E*, vol. 71, p. 036113, Mar 2005. [Online]. Available: http://link.aps.org/doi/10.1103/PhysRevE.71.036113

[35] T. Optimization, "Tomlab /cplex." [Online]. Available: http://tomopt.com/tomlab/products/cplex/

[36] X. Liu, S. Pequito, S. Kar, B. Sinopoli, and A. Aguiar, "Minimum sensor placement for robust observability of structured complex networks." [Online]. Available: https://www.mathworks.com/matlabcentral/fileexchange/57454

[37] V. Chvatal, "A greedy heuristic for the set-covering problem," *Mathematics of Operations Research*, vol. 4(3), pp. 233–235, 1979.